\documentclass{article}
\usepackage{natbib}
\usepackage{xcolor}
\usepackage{color}
\usepackage{subfigure}
\input{epsf.sty}
\usepackage{epsfig}
\usepackage{mathrsfs}  
\usepackage{amsthm,amsfonts,amssymb,amsmath}
\usepackage{multirow}
\newcommand{\NV}{firm}

\newcommand{\qeds}{$\qedsymbol $}

\newtheorem{thm}{Theorem}
\newtheorem{lem}{Lemma}
\newtheorem{pro}{Proposition}
\newtheorem{cor}{Corollary}

\def\bE{\mathbb{E}}
\def\bP{\mathbb{P}}
\def\bH{\mathscr{H}}
\def\bN{\mathscr{N}}

\begin{document}
 
\title{Cash-Flow Based Dynamic Inventory Management}
\author{Michael N. Katehakis\\
Rutgers University,  \\100 Rockafeller Road, Piscataway, NJ 08854
\and Benjamin Melamed\\
Rutgers University,  \\100 Rockafeller Road, Piscataway, NJ 08854
\and Jim (Junmin) Shi\\
New Jersey Institute of Technology\\
University Heights,   Newark, NJ 07102
}
 \date{}
\maketitle

\begin{abstract}
Small-to-medium size enterprises (SMEs), including many startup firms, need to manage interrelated flows of cash and inventories of goods.  In this paper, we model a firm that can finance its inventory (ordered or manufactured) with loans in order to meet random demand which in general may not be time stationary.  The firm earns interest on its cash on hand and pays interest on its debt.  The objective is to maximize the expected value of the firm's 
capital at the end of a finite planning horizon.  Our study shows that the optimal ordering policy is characterized by a pair of threshold variables for each period as function of the initial state of the period. 
Further,  upper and lower bounds  for the  threshold values are  developed using two simple-to-compute  ordering policies.  Based on these bounds, we provide an efficient algorithm to compute the two threshold values. 
Since the underlying state space is two-dimensional which leads to high computational complexity of the optimization algorithm, we also derive upper bounds for the optimal value function by reducing the optimization problem to one dimension.  Subsequently, it is shown that policies of similar structure are optimal when the loan and deposit interest rates are piecewise linear functions, when there is a maximal loan limit and when unsatisfied demand is backordered.  Finally, further managerial insights are provided with numerical studies.

\end{abstract}
 

\section{Introduction} \label{intro}
In the current competitive environment, small-to-medium size enterprises (SMEs), including many startup firms, must make joint decisions concerning interrelated flows of cash and products; cf. \cite{MaRS2014}.  Cash is the lifeblood of any business and is often in short supply. Cash flows stem from operations, financing and investing activities. To maximize profits, firms need to manage operations, financing and investing activities, efficiently and comprehensively.  For example, most startup firms need access to funds, typically for ordering or manufacturing products; cf. \cite{Zwilling2013}.  As another example, some small businesses are seasonal in nature, particularly retail businesses. If a business makes most of its sales during the holiday season, it may need a loan prior to the holiday season to purchase a large amount of inventory to gear up for that season. Such loans to purchase inventory are generally short-term in nature and companies usually pay them off after the season is over with the proceeds of sales from their seasonal sales.\\

To facilitate SME borrowing, \textit{Small Business Administration} (SBA) has established a loan program with a successful record.  For example, \textit{Wells Fargo Capital Finance} is the nation's leading provider (by dollar volume) of loans guaranteed by the SBA, providing loans to companies across a broad range of industries throughout the United States.  Building on a strong track record, Wells Fargo works with companies in many diverse industries, including apparel, electronics, furniture, housewares, sporting goods, toys and games, food products, hardware, and industrial goods.  Accordingly, small business loans at Wells Fargo rose 18\% in 2014; cf. \cite{Reuter2015}, \cite{WellsFargo2015} and \cite{Yahoo2015}.\\


In this paper we model and analyze the optimal financial and operational policy of an SME firm whose inventory is subject to lost sales, zero replenishment lead times and periodic review over a finite planning horizon. The firm's state (i.e., {\it inventory-capital profile}) consists of inventory level and cash level, where a positive cash level represents cash on hand while a negative cash level represents a loan position.  Cash flows are managed from the perspectives of operations, financing and investing in each period, as follows.
(i)	From the operations perspective, the firm stocks up and sells inventory;
(ii)	From the financing perspective, the firm can use its cash on hand or an external short-term loan (if needed) to procure products for inventory;
(iii)	From the investing perspective, any cash on hand is deposited in a bank account to earn interest at a given rate, while debt incurs interest at a higher given rate.
We develop and study a discrete-time model for a single-product inventory system over 
multiple periods.  The objective of the firm is to dynamically optimize the order quantities in each period, given the current state of cash and inventory, via joint operational/financial decisions, so as to maximize the expected value of the firm's capital (i.e., total wealth level) at the end of a finite time horizon.\\

The inventory flow is described as follows.  At the beginning of each period, the firm decides on an order quantity and the corresponding replenishment order materializes with zero lead time. During the remainder of the period, no inventory transactions (demand fulfillment or replenishment) take place.  Rather, all such transactions are settled at the end of the period. Incoming demand is aggregated over that period, and the total period demand draws down on-hand inventory.  However, if the demand exceeds the on-hand inventory, then all excess demand is lost (the backorder setting is studied in \S \ref{Sec:Ext:Backorder}). All the leftover inventory (if any) is carried forward to the next period subject to a holding cost, and at the last period, the remaining inventory (if any) is disposed of either at a salvage value or at a disposal cost.\\

Cash flow typically takes place as follows.  
All transactions pertaining to previous period are settled at its end.  More precisely, the firm updates its cash position with the previous period's revenue from sales and interest earned from a deposit, or paid on outstanding debt (if any).  The firm then decides the order quantity for the next period and pays for replenishment as follows: first, with cash on hand, and if insufficient, with a withdrawal from cash on deposit (there is no withdrawal penalty), and if still insufficient, by an external loan.  If not all cash on hand is used for replenishment, then any unused cash is deposited in a bank where it earns interest.  At the end of each period the resulting cash on hand or debt are carried forward to the next period.\\

The main contribution of the paper is to establish and compute the optimal ordering policy in terms of the {\it net worth} of the firm (operating capital in in product units and inventory on hand) at the beginning of each period.  It is shown that the optimal policy is characterized for each period $n$ by a pair of threshold variables, $\alpha_{n}$ and $\beta_{n}$ where $\alpha_{n}<\beta_{n}$; cf. Theorem \ref{Them:Optimal-q} for the single period problem and Theorem \ref{Them:Dynamic_order_policy} for the multi-period problem. The threshold values, $\alpha_{n}$ and $\beta_{n}$, are in general functions of the firm's net worth. This optimal policy has the following structure:
a)	If the net worth is less than $\alpha_n$, then the firm orders up to $\alpha_n$, a case referred to as the {\it over-utilization} case.
b)	If the net worth is greater than $\beta_n$, then the firm orders up to $\beta_n$, a case referred to as the {\it under-utilization} case.
c)	Otherwise, when the net worth is between $\alpha_{n}$ and $\beta_{n}$, then the firm orders exactly as many units as it can afford {\it without borrowing}, a case referred to as the {\it full-utilization} case.
%
For the single period problem, we further derive the optimal solution in closed form and show that the $(\alpha,\beta)$ optimal policy yields a positive expected value even with zero values for both initial inventory and cash. For the multi-period problem, we construct two myopic policies, which provide upper and lower bounds for the threshold values, respectively.  Based on these upper and lower bounds, we provide an algorithm for computing the two thresholds, $\alpha_n$ and $\beta_n$, recursively for all periods $n$. 
Since the underlying state space is two-dimensional, it leads to high computational complexity of the optimization algorithm, we   derive upper bounds for the optimal value function by reducing the optimization problem to one dimension.
Subsequently, it is shown that policies of similar structure are optimal when the loan and deposit interest rates are piecewise linear functions, when there is a maximal loan limit and when unsatisfied demand is backordered; cf. Theorems \ref{Corr:Optimal-q-piecewise}, \ref{Them:Dynamic_order_policy_Loan_Capacity} and \ref{Them:Dynamic_order_policy-backorder}, respectively.\\

The remainder of this paper is organized as follows. Section \ref{literature} reviews related literature and Section \ref{Sec:N-period-model} formulates the model. In Section \ref{Sec:single_models}, the single-period model is developed and the optimal ($\alpha,\beta$) policy is derived, where the threshold values are functions of the demand distribution and the cost parameters of the problem. Section \ref{Sec:N-period-model-analysis} extends the analysis for the dynamic multi-period problem and derives a two-threshold structure, ($\alpha_{n},\beta_{n}$) policy, for the optimal policy via a dynamic programming analysis. Section \ref{Sec:myopic} introduces two myopic policies that provide upper and lower bounds for each $\alpha_n$ and $\beta_n$, based on which an efficient searching algorithm for $\alpha_{n}$ and $\beta_{n}$ is presented. An upper bound for the value function is also introduced. Numerical studies are presented in Section \ref{Sec:numerical study}.  In Section \ref{Sec:discussion}, it is pointed out that simple modifications of the underlying model allow the extension of the results to the case of piecewise linear interest rate functions, a maximal loan limit and backorder of unsatisfied demand. Finally, Section \ref{Sec:conclusions} concludes the paper. 
All proofs are relegated to the appendix.

\section{Literature Review} \label{literature}

In a seminal study, \cite{MM1958} show that in a perfect capital market with adequate and inexpensive external funding, a firm's operational and financial decisions can be made separately. Since then, most literature in inventory management extended the classical newsvendor problem in a variety of ways, but assuming that decision makers are not subject to financial constraints. However, in view of the imperfection of real-life capital market, a growing body of literature has begun to consider operational decision making subject to financial constraints. In those studies, inventories of goods are often treated as special financial instruments; cf. \cite{Singhal1988}. Accordingly, portfolios of physical products and financial instruments have been studied using finance/investment principles such as {\it Modern Portfolio Theory} (MPT) and {\it Capital Asset Pricing Model} (CAPM). For related literature, the reader is referred to \cite{Corbett1999} and references therein. Further, the relationship between inventories and finance, along with their theoretical and empirical consequences, are discussed in \cite{Girlich2003}.\\

The growing literature on the interface between operations management and financial decisions can be categorized into two major strands: \textit{single-agent} and \textit{game related multiple-agent} streams. Each stream, attendant models can be further classified as single-period and multi-period models. Our study belongs to the single-agent strand and addresses both single-period and multi-period models.\\

In the single-agent strand, a firm operational decisions and financial decisions are typically made simultaneously, but without interacting with other firms.  Considering an imperfect market, \cite{XuBirge2004} develop models with simultaneous production and financing decisions in the presence of demand uncertainty, which illustrates how a firm's production decisions are affected by the existence of financial constraints. Recently, \cite{Birge2011} present an extension of a model in \cite{XuBirge2004} by assuming that debt and production scale decisions have fixed costs necessary to maintain operations, variable costs of production, and volatility in future demand forecasts.  Building on their previous work of \cite{XuBirge2004}, \cite{XuBirge2008} review that analysis and consider the effect of different operating conditions on capital structure, and exhibit some empirical support of their previously predicted relationship between production margin and market leverage.  \cite{Cai2014Credit} investigate the roles of bank and trade credits in a supply chain with a capital-constrained retailer facing demand uncertainty. The studies above focus on single-period problems.\\

Some single-agent, multi-period problems are addressed in \cite{Hu2005capital} and \cite{Li2013control}. Here, multi-period models posit different interest rates on cash on hand and outstanding loans. These papers also demonstrate the importance of the joint consideration of production and financing decisions in a start-up setting in which the ability to grow the firm is mainly constrained by its limited capital and dependence on bank financing. For example, \cite{Hu2005capital} examine the interdependence of a firm's capital structure and its short-term operating decisions concerning inventories, dividends, and liquidity. To this end, \cite{Hu2005capital} formulate a dynamic model to maximize the expected present value of dividends. \cite{Li2013control} study a dynamic model of managerial decisions in a manufacturing firm in which inventory and financial decisions interact and are coordinated in the presence of demand uncertainty, financial constraints, and default risk. It is shown that the relative financial value of coordination can be made arbitrarily large.  In addition, \cite{Gupta2009tradecredit} study a retailer's dynamic inventory problem (both discrete and continuous review) in the presence of random demand while the retailer is financially supported by trade credit from its supplier, contingent on the age of the inventory. It is shown that the optimal policy is a base stock policy where the base stock level is affected by the offered credit terms.  \cite{Luo2013Cash} study a centralized supply chain consisting of two divisions (modeled as echelons), of which the headquarter division manages the financial and operational decisions with a cash pool. They focus on analyzing the value of cash pooling.\\

The work most closely related to this study is \cite{chao2008dynamic} and \cite{Gong2014dynamic}. The first reference considers a single-agent multi-period problem for a self-financing retailer without external loan availability. It shows that the optimal, cash flow-dependent, policy in each period, is uniquely determined by a single critical value.  Our study differs from \cite{chao2008dynamic} in several ways.  First, from a modeling perspective, \cite{chao2008dynamic} consider a self-financing firm without access to external credit, while we consider possible loan which provides the retailer with the flexibility of ordering a larger quantity to achieve a higher fill rate. The consideration of loan intertwines the financial decision (e.g., how much to borrow) with the operational decision (e.g., how many to order) as well as the investment decision (how much to deposit for additional interest). In this case, the complexity of the problem is raised from one dimension (operational decision alone) to two dimensions (joint operational and financial decisions).  Further, \cite{chao2008dynamic} assume an iid demand process and time-stationary costs, while we consider a non-stationary demand process and time varying (loan and deposit) interest rates.  Second, in terms of results, \cite{chao2008dynamic} structure the optimal policy as a base stock policy with a single threshold value; while our paper shows that the optimal policy is characterized by two threshold values, which divide the ordering decision space into three intervals to guide financing, operational and investment decisions.  
\cite{Gong2014dynamic} study a similar  model to maximize the expected terminal total capital but along different lines from ours.
 However, this paper differs in several aspects as follows.  In terms of model setting, we consider non-stationary demand process and time varying (loan and deposit) interest rates as well as an extension to realistic piecewise loan rates.  In terms of results, our optimal policy is presented in terms of two thresholds per period (multiple thresholds for the case of  piecewise loan rates), while \cite{Gong2014dynamic} analyze  a single-threshold policy.  Computationally, we develop simple easy-to-compute policies (based on single period information) which provide lower and upper bounds for each of the two threshold values, and we propose an efficient algorithm for the computation of the optimal policy. Further we provide upper and lower bounds for the value function.  Finally we develop further managerial insights into the interplay between operational, financial and investment decisions.\\

In a second strand of the literature, multi-agent competition between firms and financial institutions has been investigated using game theoretic approaches. This literature includes, but is not limited to, \cite{buzacott2004inventory}, \cite{dada2008},  \cite{Yasin2010} and \cite{Raghavan2011lending}. Most such papers deal with a single-period problem.  \cite{buzacott2004inventory} analyzes a Stackelberg game between the bank and the retailer in a newsvendor inventory model.  \cite{dada2008} assume that the interest rate is charged by the bank endogenously and use a game model to capture the relation between the bank and the inventory controller through which an equilibrium is derived and a non-linear loan schedule is obtained to coordinate the channel.  \cite{Yasin2010} study the implications of asset based lending for operational investment, probability of bankruptcy, and capital structure for a borrower firm.  \cite{Raghavan2011lending} study a short-term financing problem in a cash-constrained supply chain.\\

As a part of the second literature strand, multi-agent game-theoretic approaches have also been used to model the competition between suppliers and retailers.  Such recent studies include \cite{Kouvelis2011a} and \cite{Kouvelis2011b} among many others.  An important portion of this strand addresses the impact of trade credits provided by suppliers to retailers.  \cite{lee2010fiancingcost} study the impact of inventory financing costs on supply chain coordination by considering four coordination mechanisms: all-unit quantity discount, buy backs, two-part tariff, and revenue-sharing. It is shown that using trade credits in addition to contracts, a supplier can fully coordinate the supply chain and achieve maximal joint profit. Further, \cite{lee2011tradecredit} model a firm with a supplier that grants trade credits and markdown allowances. Given the supplier's offer, it determines the order quantity and the financing option for the inventory under either trade credit or direct financing from a financial institution.   The impact of trade credits from an operational perspective is also studied in \cite{Yang2011inventory}, which investigates the role that trade credit plays in channel coordination and inventory financing.  It is shown that when offering trade credits, the supplier balances its impact on operational profit and costs. 
Other related work of ours includes:
\cite{shi2014production}, \cite{katehakis2012computing}, 
\cite{zhou2007effective}, \cite{zhao2006structure}.

\section{Model Formulation} \label{Sec:N-period-model}
Consider a time horizon $\bN:=\{1,2,\ldots, N\}$ where the periods are indexed forward.
At the beginning of period $n\in \bN$,
let the ``\textit{inventory-capital}''  state of the system
be characterized by a vector $(x_n, y_n),$
where $x_n\geq 0$ denotes the amount of on-hand inventory (number of product units) and
$y_n$ denotes the amount of product that can be purchased using  all the available capital, i.e., $y_n$ is the capital position measured in ``product units''.
Let $q_n\geq 0$ denote the order  quantity the \NV\ uses in the beginning of period $n\in \bN .$ The ordering cost per unit is $c_n\geq 0$.
We assume that any order quantity can be fully satisfied within zero lead time of replenishment.
Let $D_n\in \Re^{+}$ denote the random demand during period $n$ which follows a general distribution.
We assume that demands of different
periods are independent but could be generally non-stationary across periods.
Let  $f_n(\cdot)$ and $F_n(\cdot)$ denote the \textit{probability density function} (pdf) and  the \textit{cumulative distribution function} (cdf) of $D_n,$ respectively. The selling price per unit is $p_n\geq c_n$.\\

Throughout all
periods $n=1,\ldots,N-1$, any unsold units are
carried over in inventory to be used in subsequent periods subject
to a constant holding cost $h_n\geq 0$ per unit over period $n$,  and any unmet demand is lost. The backorder setting of the problem is studied in \S \ref{Sec:Ext:Backorder}.
%
At the end of the horizon, i.e., the end of period $N$, all leftover
inventory (if any) will be salvaged (or disposed off) at a constant price
(cost) $s$ per unit. Note that we allow a negative $s$ in which case $s$ represents a disposal cost per unit,  e.g., the unit cost of disposing vehicle tires. For notational convenience, we denote $h_N=-s$ in the sequel, so that the salvage value (or disposal cost) $s$ is treated uniformly as a type of inventory cost.\\

To finance inventory, a loan can be applied at interest rate $\ell_n\geq 0$. Interest rates are fixed and depend entirely upon current market conditions. 
  For example, \textit{Lending Club} provides Small Business Loans at a fixed rate financing up to \$300,000
 (Lending Club is a US peer-to-peer lending company, headquartered in San Francisco, California).
When there remains unused capital after an ordering decision, the unused capital is deposited in a bank account to earn interest at rate $i_n\geq 0$.
To avoid trivialities, we assume that $i_n<\ell_n$ and it is
possible to achieve a positive profit with the aid of a loan, i.e., $(1+\ell_n)c_n < p_n$.
This assumption can be equivalently written as
$1+\ell_n<p_n/c_n$.
Note also, that the above assumption implies $1+i_n < p_n/c_n$ since $i_n<\ell_n$;
namely, investing on inventory is more profitable than depositing available capital in
the bank.\\
%


Given the initial inventory-capital state as $(x_n,y_n),$ if an order of size
 $q_n \ge 0 $ is placed and the demand during the period is $D_n$, then we have the following operational and financial cash inflows:
 
\begin{itemize}
 \item
The cash flow from operational sales of items (i.e., the
realized revenue from  inventory) at the end of the period is given by
\begin{eqnarray}
R_n(D_n,q_n,x_n)&=&p_n\cdot \min\{q_n+x_n,D_n\}-h_n\cdot [q_n+x_n-D_n]^{+} \nonumber\\
  &=&p_n\cdot[q_n+x_n-(q_n+x_n-D_n)^+]-h_n\cdot [q_n+x_n-D_n]^{+} \nonumber\\
    &=&p_n(q_n+x_n)-(p_n+h_n)\cdot [q_n+x_n-D_n]^{+},\label{Eq:Inv_Revenue}
\end{eqnarray}

 \noindent where $[a]^+=\max\{a,\,0\}$ denotes the positive part of
real number $a$, and the second equality holds by
$\min\{a,b\}=a-[a-b]^+$.
\\
\item
The cash flow from capital (i.e., the realized revenue from financial investment) at the end of the period
is computed for each of the following investing and financing scenarios:
\begin{itemize}
 \item[i)]
Investing decision: If the order quantity satisfies $0 \leq q_n\leq y_n$, then the leftover amount $c_n\cdot(y_n- q_n)$ of cash is deposited in the bank and will yield a positive inflow of $c_n \cdot (y_n-q_n)(1+i_n)$  at the end of the period $n$.

\item[ii)] Financing decision:  If the order quantity satisfies $q_n > y_n $ (including the case $q_n=0>y_n$), then a loan amount of $c_n\cdot (q_n-y_n)$ will be borrowed during the period and will result in a cash outflow of $c_n \cdot (q_n-y_n)(1+\ell_n)$ at the end of the period.

\end{itemize}
In summary, the cash flow from the bank  (positive or negative)  can
be written in general as
\begin{eqnarray}
K_n(q_n,y_n)  =c_n \cdot  (y_n-q_n)\left[(1+i_n){\bf 1}_{\{q_n\leq y_n\}}+(1+\ell_n){\bf
1}_{\{q_n> y_n\}}\right].
      \label{Eq:MarketRevenue}
\end{eqnarray}
\end{itemize}

Note that the cash inflow  from inventory, $R_n(D_n,q_n,x_n) $  is dependent on $D_n$ while
independent of  $y_n$. In a similar vein, the cash flow from capital, $K_n(q_n,y_n)$  is independent of the
initial on-hand inventory, $x_n$ and the demand size $D_n$, but is dependent on $y_n$.
Also, note that the ordering cost, $c_n\cdot q _n$, has been accounted for in Eq. (\ref{Eq:MarketRevenue}) while the remaining capital, if any,  is invested in the bank and its value at the end of the
period is given by  $K_n(q_n,y_n)$.\\

Since the order quantity $q_n=q_n(x_n,y_n)$ is decided at the
beginning of period $n$ as a function of $(x_n,y_n)$, it is readily
shown that the state process $\{(x_n,y_n)\}$ under study is a {\it
Markov decision process} (MDP) with decision variable $q_n$; cf.
Ross (1992). 
The state dynamics of the system, i.e., inventory flow and cash flow,  are given as follows, for $n=1,2,...,N-1$
\begin{eqnarray}
x_{n+1}&=&[x_{n}+q_{n}-D_{n}]^{+}; \label{Eq:xn0}\\
y_{n+1}&=&
[R_{n}(D_n, q_{n},x_{n})+K_{n}(q_{n},y_{n})]/c_{n+1},
\label{Eq:yn0}
\end{eqnarray}
where $R_n$ and $K_n$ refer to the capital gain from operational and financial decisions, as given by Eqs. (\ref{Eq:Inv_Revenue})-(\ref{Eq:MarketRevenue}), respectively.\\
 

Note that, at the beginning of period $n$, the interplay between finance and operation is one way. Namely, it is feasible to purchase products
with available capital $y_n$ (when  $y_n > 0$),  but it is typically restricted to convert any on-hand inventory $x_n$  into cash. Conversely, the sales interplay between finance and operations is opposite, namely, sales deplete the inventory and contributes to cash on hand.  

Thus, we have the following dynamic programming formulation:
\begin{eqnarray}
V_n(x_{n},y_{n})= \sup_{q_n \geq 0}{\bE}\bigg[V_{n+1}(x_{n+1},y_{n+1}) |x_{n},y_{n}\bigg], \hspace{20 pt}
n=1,2,\cdots,N-1 \label{Eq:DP-optimal}
\end{eqnarray}
where  the expectation is taken with respect to $D_n$, and $x_{n+1} ,$  $y_{n+1}$ are given by Eqs. (\ref{Eq:xn0}),
(\ref{Eq:yn0}), respectively. For the final period $N,$ we have
\begin{eqnarray}\label{Eq:DP-optimal-N}
V_N(x_{N},y_{N})= \sup_{q_N \geq 0}{\bE}\bigg[
R_{N}(D_N,q_{N},x_{N})+K_{N}(q_{N},y_{N}) \bigg].
\end{eqnarray}

\section{The Single Period Problem}
\label{Sec:single_models}

%
The analysis of the last period $N$ is a typical single-period problem and its objective function is given by Eq. (\ref{Eq:DP-optimal-N}).
For notational convention, we omit the subscript $N$ for each related variable when denoting counterparts under the single period setting, without causing any confusion.

Accordingly, for any given initial ``inventory-capital''  state $(x,y)$ and ordering decision $q\geq 0$, the expected value of net worth at the end of the period is given by
%
%
%
\begin{eqnarray}
G(q,x,y)=p\,(x+q)-(p-s) \int_0^{x+q}(x+q-t)f(t)dt 
+c \, (y-q)\bigg[(1+i){\bf 1}_{\{q\leq y\}}+(1+\ell){\bf 1}_{\{q>
y\}}\bigg]. \label{Eq: E profit}
\end{eqnarray}

The following lemma summarizes the important properties of the function $G(q,x,y)$.
\begin{lem}\label{Lem:concave}
Function $G(q,x,y)$ is continuous in $q,$ $x$ and $y ,$ and
it has the following properties. 
\begin{enumerate}
 \item[i)] It  is  concave in $q\in
[0,\infty)$, for all $x ,$   $y$ and  all $s <p  .$
\item[ii)] It is increasing and concave  in $x,$ for $s \geq 0. $
\item[iii)] It is increasing and concave  in $y,$ for all $s <p.$
\end{enumerate}
\end{lem}

\bigskip

{\bf Remarks.}
\begin{enumerate}
 \item
It is important to point out that $G(q,x,y)$ might not increase in
$x$ if $s<0 .$ In particular,   if $s$ represents a disposal cost,
i.e., $s<0$, then 
$\partial G(q,x,y)/\partial x$ might
be negative, which implies that $G(q,x,y)$ is decreasing for some
high values of $x . $

Further, for the special case $s<0,$ it is of interest to locate the
critical value, $x'$ such that $G(q,x,y)$ is decreasing for $x>x'$.
To this end, we set Eq. (\ref{Eq: derv profit_x}) to be zero, which
yields
\begin{eqnarray}
(p-s)\cdot F(q+x)=p \label{Eq: derv profit_x=0} .
\end{eqnarray}
Therefore,
\begin{eqnarray}
x'=F^{-1}\left(\frac{p}{p-s}\right)-q \label{Eq: derv profit_x'},
\end{eqnarray}

where $F^{-1}(\cdot)$ is the inverse function of
$F(\cdot)$. Eq. (\ref{Eq: derv profit_x'}) shows that a higher disposal cost,
$-s$, implies a lower threshold for $x'$ above.

\item
 Lemma \ref{Lem:concave} implies
that  higher values of initial assets, $x$ and $y$ or the net  worth $\xi=x+y$, will
yield a higher expected revenue $G(q,x,y) $. Further, for any
assets $(x,y)$ there is a unique optimal order quantity $q^{*}$ such
that
$q^{*}(x,y)={\arg
\max}_{q \geq 0}G(q,x,y).$
\end{enumerate}

We next define  the critical values, $\alpha$ and
$\beta$, as follows:
\begin{eqnarray}
\alpha&=&F^{-1}(a) ; \label{Eq:alpha}\\
\beta&=&F^{-1}(b) , \label{Eq:beta}
\end{eqnarray}
\noindent where
$$a=\frac{p-c\cdot (1+\ell)}{p-s} ; \,\,
b=\frac{p-c\cdot (1+i)}{p-s} .$$
%
 
It is readily seen that  $a\leq b ,$ since $0\leq i \leq \ell$ by assumption. This implies that  $\alpha\leq \beta ,$
since $F^{-1}(z)$ is increasing in $z$. The critical value $\beta$ can be interpreted as the optimal
order quantity for the classical \NV\ problem corresponding to the case of sufficiently large
 $Y$ of our model, in which case no loan is involved, but the unit ``price'' $c\cdot (1+i)$ has been inflated to reflect the opportunity  cost
of cash not invested in the bank at an interest rate $i.$  Similarly,  $\alpha$ can
be interpreted as the optimal order quantity for the classical \NV\
problem corresponding to the case $Y=0 $ of our model, i.e., all units are purchased by a loan at an interest rate $\ell .$\\

Note also that in contrast to the classical \NV\ model, the critical values
 $\alpha$ and $ \beta $ above,  are now functions of the corresponding interest rates and represent opportunity costs that take into account the value of  time using the interest factors $1+i$ and $1+ \ell$.\\

We next provide the following theorem regarding the
optimality of the  ($\alpha$, $\beta$)  policy.
 
\begin{thm}\label{Them:Optimal-q}
For any given initial inventory-capital state $(x,y)$, the optimal order quantity is
\begin{eqnarray}\label{eqstar}
q^{*}(x,y)= \left\{ \begin{array}{ll}
            (\beta-x)^+, &\mbox{  $ \beta \leq x+y$}; \\[0.3cm]
            y^+, &\mbox{  $\alpha\leq x+ y < \beta$}; \\[0.3cm]
  (\alpha-x)^+, &\mbox{  $x+y < \alpha$,}
       \end{array} \right.
\end{eqnarray}
where $\alpha$ and $\beta$ are given by Eq. (\ref{Eq:alpha}) and
(\ref{Eq:beta}), respectively.
\end{thm}
%
 
Note that the optimal ordering quantity in the classical newsvendor model [cf.
\cite{zipkin2000foundations} and many others], can be obtained from Theorem \ref{Them:Optimal-q}  as the solution for the extreme  case with  $i=\ell=0,$ in which case the optimal order
quantity is given by:
$$\alpha=\beta=F^{-1}\left(\frac{p-c}{p-s}\right). $$

\begin{figure}{}
  \centering
   \includegraphics[width=0.5\textwidth]{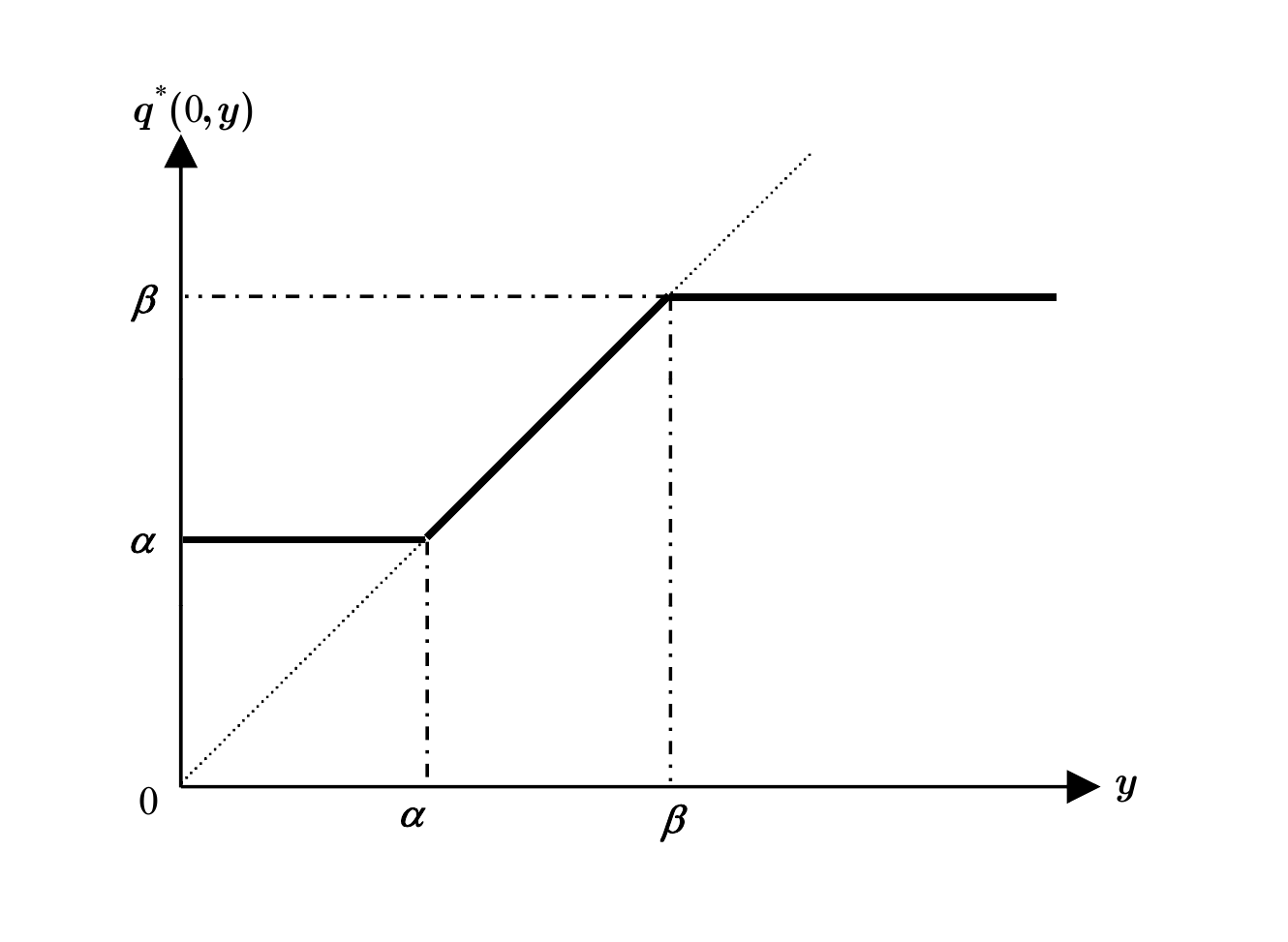}
\caption{The Optimal
Order Quantity when $x=0$} \label{Fig:Q*=0}
\end{figure}

Next we elucidate the structure of the ($\alpha,\beta$) optimal policy below where we discuss the utilization level of the  initially available capital $Y$.
\begin{enumerate}
\item {\bf Over-utilization}: When $x+y<\alpha$, it is optimal to order $q^*=\alpha-x = y +(\alpha -x -y) $. In this case  $y=Y/c$ units are bought using
 all the available fund $Y$ and the remaining $(\alpha -x -y)$ units are bought using a loan of size: $c \cdot (\alpha -x -y) .$
\item {\bf Full-utilization}: When $\alpha\leq x+y<\beta$, it is optimal to order $q^*=y =Y/c$ with all the
available fund $Y$. In this case, no deposit and no loan are involved.
\item {\bf Under-utilization}: When $x+y\geq \beta$, it is optimal
to order $q^*=(\beta-x)^+$. In particular, if  $x< \beta$, it is optimal to order
$\beta-x$   using $c \cdot (\beta-x)$ units of the available cash $Y$, and deposit the remaining cash to earn interest. However, if $x \geq \beta$, then  $q^*=0 $, i.e., it is optimal not to order any units and deposit all the amount of $Y$ to earn interest.
\end{enumerate}
The above interpretation is illustrated in Figure \ref{Fig:Q*=0}  for the case  $x=0, $ which plots  the optimal order quantity $q^*$ as a function of $y$. Note that for $y \in (0,\alpha)$ there is over utilization of  $y$; for $y \in [\alpha,\beta)$ there is full utilization of $y$ and for
$y \in [\beta, \infty)$ there is under utilization of  $y$.

In light of aforementioned $(\alpha,\beta)$ policy, one can derive the optimal expected net worth, $V(x,y)= \max_{q\ge0} G(q,x,y)$  in explicit form.\\

\begin{thm}\label{Them:Optimal-v} For the single-period problem with initial state $(x,y)$,

i) V(x,y) is given by
\begin{eqnarray}\label{Eq:V_star}
V(x,y)= \left\{ \begin{array}{ll}
            p\cdot x-(p-s)\cdot T(x)+c\cdot y\cdot (1+i), &\mbox{  $x>\beta$}; \\[0.3cm]
            p\cdot \beta-(p-s)\cdot T(\beta)+c\cdot (x+y-\beta)(1+i), &\mbox{ $x\leq \beta$,  $ \beta \leq x+y$}; \\[0.3cm]
            p\cdot (x+y^+)-(p-s)\cdot T(x+y^+), &\mbox{  $\alpha\leq x+ y < \beta$}; \\[0.3cm]
  p\cdot \alpha-(p-s)\cdot T(\alpha)+c\cdot (x+y-\alpha)(1+l), &\mbox{  $x+y < \alpha$,}
       \end{array} \right.
\end{eqnarray}
\indent where $T(x)=\int_{0}^{x}(x-t)f(t)dt$;

ii) The function $V(x,y)$
is increasing in $x$ and $y$, and jointly concave in $(x,y)$, for $x,y\geq 0$.
\end{thm}

\bigskip
From an investment perspective, it is of interest to consider the possibility of speculation; cf. \cite{Hull2002}. The following result shows that the $(\alpha,\beta)$ policy given in Theorem \ref{Them:Optimal-q} yields positive value with zero investment. Specifically, when the
\NV\ has zero initial inventory and capital, i.e., $x=0$ and $y=0$, the optimal policy yields a positive expected final asset value.

\begin{cor}\label{Cor:speculation} For $x=0$ and $y=0$, one has
   $$V(0,0)=(p-s)\int_0^{\alpha}t\cdot f(t)dt>0.$$
\end{cor}

Note that arbitrage usually means that it is possible to have a positive profit
for any realized demand (i.e., of a risk-free profit at zero cost; cf. \cite{Hull2002}), thus the above speculation possibility does not in general imply an arbitrage, which only exists if the demand is constant.

\section{Optimal Solution for the Multiple Period Problem}
\label{Sec:N-period-model-analysis}


For ease of exposition, we shall denote $\xi_n=x_n+y_n$ and take
$z_n=x_n+q_n$ as the decision variable, in lieu of $q_n$.  Here, $z_n$ refers to the available inventory after replenishment, and it is constrained by $z_n\geq x_n\geq 0$ for each period $n$. In this fashion, the MDP model defined by Eqs. (\ref{Eq:DP-optimal})-(\ref{Eq:DP-optimal-N}) can be presented as:
\begin{equation}
V_n(x_{n},y_{n})=\max_{z_n\geq x_n}G_n(z_n,\xi_{n}),\label{Eq:G-V} \\
\end{equation}
where $G_n(z_n,\xi_{n})={\bE}\left[V_{n+1}(x_{n+1},y_{n+1})\right]$.
The inventory-capital states are dynamically updated as
\begin{eqnarray}
x_{n+1}&=&[z_{n}-D_{n}]^{+};  \label{Eq:xn}\\
y_{n+1}&=&
p'_{n}\cdot z_{n}-(p'_n+h'_n)
\left[z_{n}-D_{n}\right]^{+}
+c'_{n} \cdot
(\xi_{n}-z_{n})\bigg[(1+i_n){\bf 1}_{\{z_{n}\leq \xi_n\}}+(1+\ell_n){\bf
1}_{\{z_{n}> \xi_n\}}\bigg],
\label{Eq:yn}
\end{eqnarray}
where
$p'_n=p_n/c_{n+1}$, $h'_n=h_n/c_{n+1}$ and
$c'_n=c_n/c_{n+1}$. %
%

We first present the following result.
\begin{lem}\label{Lem:increasing}
For $n \in \bN$, the following is true:\\
(1)  $G_n(z_n,x_{n},y_{n})$ is increasing in $x_n$ and $ y_{n}$, and is
concave in $z_n$, $x_{n}$ and
$y_{n}$.\\
(2) $V_n(x_{n},y_{n})$ is increasing and concave in $x_{n}$ and $y_{n}$.
\end{lem}

\bigskip
We next present and prove the main result of this section.

\begin{thm}\label{Them:Dynamic_order_policy} (\textbf{The ($\alpha_n$, $\beta_n$) ordering
policy}).\\
For every period $n\in \bN$ with state $(x_n,y_n)$ at the
beginning of the period, there exist positive constants
$\alpha_n=\alpha_n(x_n,y_n)$ and $\beta_n=\beta_n(x_n,y_n)$ with
$\alpha_n\leq \beta_n$, which give rise to the  optimal order quantity  as
follows:
\begin{eqnarray}\label{Eq:q*}
q_n^{*}(x_n,y_n)= \left\{ \begin{array}{ll}
            (\beta_n-x_n)^+, &\mbox{  $\xi_n\geq \beta_n$}; \\[0.3cm]
             (y_n)^+, &\mbox{  $\alpha_n\leq \xi_n < \beta_n$}; \\[0.3cm]
             (\alpha_n-x_n)^+, &\mbox{  $\xi_n < \alpha_n$.}
       \end{array} \right.
\end{eqnarray}
Further,  ${\alpha}_n$ is uniquely identified by
\begin{eqnarray}\label{Eq: N_derv profit_qi}
\bE \left[\left( \frac{\partial V_{n+1}}{\partial x_{n+1}}-(
p'_n+h'_n)\frac{\partial V_{n+1}}{\partial y_{n+1}}\right) {\bf
1}_{\{\alpha_n>D_n\}}\right]=\left[c'_n(1+\ell_n)-p'_n\right]{\bE}\left[\frac{\partial V_{n+1}}{\partial y_{n+1}}\right] ,
\end{eqnarray}
and ${\beta}_n$ is uniquely identified by
\begin{eqnarray}\label{Eq: N_derv profit_ql}
\bE \left[\left( \frac{\partial V_{n+1}}{\partial x_{n+1}}-(
p'_n+h'_n)\frac{\partial V_{n+1}}{\partial y_{n+1}}\right) {\bf
1}_{\{\beta_n>D_n\}}\right]=\left[c'_n(1+i_n)-p'_n\right]{\bE}\left[\frac{\partial V_{n+1}}{\partial y_{n+1}}\right] .
\end{eqnarray}
where the expectations are taken with respect to $D_n$ conditionally on the initial state $(x_n,y_n)$.
\end{thm}
\bigskip

Theorem \ref{Them:Dynamic_order_policy}  establishes  that the optimal ordering policy is characterized by two threshold variables. More importantly, these two threshold values, $\alpha_{n}$ and $\beta_{n}$, can be obtained recursively by solving  the implicit equations, Eqs. (\ref{Eq: N_derv profit_qi}) and  (\ref{Eq: N_derv profit_ql}), respectively.  

\bigskip
\noindent {\bf Remark.} The study in \cite{chao2008dynamic} assumes that borrowing is not allowed, and thus the \NV\ is firmly limited to order at most $y_n$ units for period $n$.  In that paper, it was shown that the optimal policy is determined in each period by one critical value. Our results presented in Theorem \ref{Them:Dynamic_order_policy} subsume that  study as a special case. This can be seen if we set $\ell_n$ to be sufficiently large such that a loan is financially prohibited. In this case, $\alpha_n$ becomes zero and $\beta_n$ can be interpreted as the critical value developed in \cite{chao2008dynamic}.
\begin{cor}\label{Cor:alpha-beta-x+y}
For any period $n<N$ and its initial state $(x_n,y_n)$, 
the threshold variables of $\alpha_n$ and $\beta_n$ are only determined by the total worth $\xi_n = x_n+y_n$, i.e., they are of the form: $\alpha_n=\alpha_n(\xi_n)$ and $\beta_n=\beta_n(\xi_n)$. But for the last period $N$, $\alpha_N$ and $\beta_N$ are independent of either $x_N$ or $y_N$.
\end{cor}
Not surprisingly, Corollary \ref{Cor:alpha-beta-x+y} implies that the operational decisions are affected by the firm's current financial status.
In view of Corollary \ref{Cor:alpha-beta-x+y}, one can first compute $\alpha_n$ and $\beta_n$ at the beginning of the period based on the total worth $\xi_n$. The decision on the order quantity $q_n$ can then be made via Eq. (\ref{Eq:q*}).

\section{Myopic Policies and Bounds Analysis}
\label{Sec:myopic}
As shown by Theorem \ref{Them:Dynamic_order_policy},   the computation of $\alpha_n$ and $\beta_n$ is complex and costly.
This fact motivates the study to follow myopic ordering policies that are relatively simple to implement. 
Such myopic policies optimize a given objective function with respect to any single period and ignore multi-period interactions and cumulative effects.  To this end, we introduce two types of myopic policies. Specifically, myopic policy (I)
 assumes the associated cost for the leftover inventory $\underline{s}_{n}$ is only the holding cost,
i.e., $\underline{s}_{n}=  -h_n.$
Myopic policy (II) assumes that  the leftover inventory  cost  $\bar{s}_{n}$ is not only the holding cost but it also includes its value in the following period,
i.e., $\underline{s}_{n}=  c_{n+1}-h_n $.
In the following two subsections, we will show that myopic policy (I) yields lower bounds, $\underline \alpha_n$ and $\underline \beta_n$ for the two threshold values, $\alpha_n$ and $\beta_n$, while myopic policy (II) yields upper bounds, $\bar \alpha_n$ and $\bar \beta_n$.\\

Before presenting the myopic policies, we present the following lemma that will be applied to derive the upper and lower bounds. 
\begin{lem}\label{Lem:E-inequality}
For real functions $f(x)$ and $g(x)$,\\
\noindent (a) if both $f(x)$ and $g(x)$ are monotonically increasing or decreasing, then
$$\bE[f(X)\cdot g(X)]\geq \bE[f(X)]\cdot \bE[g(X)],$$
where the expectation is taken with respect to the random variable  $X$.\\
\noindent (b) If $f(x)$ is increasing (decreasing), while $g(x)$ is decreasing (increasing), then
$$\bE[f(X)\cdot g(X)]\leq \bE[f(X)]\cdot \bE[g(X)].$$
\end{lem}

 \subsection{Myopic Policy (I) and Lower      Bounds for $\alpha_n$ and $\beta_n$}
\label{Sec: Lower Bounds}

Myopic policy (I) is the one period optimal policy obtained when we change
the periodic cost structure by assuming that only the  holding cost is assessed for any
leftover  inventory i.e., we assume the following modified ``salvage value'' cost structure:
\begin{eqnarray}\nonumber
{\underline s}_n = \left\{ \begin{array}{ll}
             -h_n, &\mbox{}n<N , \\
             s, &\mbox{}n=N.
       \end{array} \right.
\end{eqnarray}
%
\noindent Let further,
\begin{eqnarray}
{\underline a}_n&=&\frac{p_n-c_n \cdot (1+\ell_n)}{p_n-{\underline s}_n}\label{Eq:alpha_hat};\\
{\underline b}_n&=&\frac{p_n-c_n \cdot (1+i_n)}{p_n-{\underline s}_n}. \label{Eq:beta_hat}
\end{eqnarray}
and the corresponding critical values are respectively given by
\begin{eqnarray}
{\underline \alpha}_n&=&F^{-1}_{n}({\underline a}_n) \label{Eq:myoptical-alpha-hat}; \\
{\underline \beta}_n&=&F^{-1}_{n}({\underline b}_n). \label{Eq:myoptical-beta-hat}
\end{eqnarray}
\noindent  For $n\in \bN,$ the   order quantity below defines   the {\it myopic
 policy} (I):
\begin{eqnarray}\nonumber
{\underline q}_n(x_n,y_n)= \left\{ \begin{array}{ll}
            ({\underline \beta}_n-x_n)^+, &\mbox{  $x_n+y_n\geq {\underline \beta}_n$}; \\[0.3cm]
             (y_n)^+, &\mbox{  ${\underline \alpha}_n\leq x_n+ y_n < {\underline \beta}_n$}; \\[0.3cm]
             {\underline \alpha}_n-x_n, &\mbox{  $x_n+y_n < {\underline \alpha}_n$.}
       \end{array} \right.
\end{eqnarray}
The next theorem establishes the lower bound properties of the myopic policy (I).
\bigskip

\begin{thm}\label{Thm:Optimal-myopical-hat}
The following are true:\\
i) For the last period $N$, $\alpha_N={\underline \alpha}_N$ and $\beta_N = {\underline \beta}_N$.\\
ii) For any period $n=1,2,\ldots N-1$,
$\alpha_n  \geq  {\underline \alpha}_n$ and
$\beta_n   \geq  {\underline \beta}_n .$
%
%
\end{thm}

\subsection{Myopic Policy (II) and Upper Bounds for $\alpha_n$ and $\beta_n$}

\label{Sec: Upper Bounds}


Myopic policy (II) is the one period optimal policy obtained when we change
the periodic cost structure by assuming that not only the  holding cost is assessed but also the cost in the next period for any leftover inventory i.e., we assume the following modified ``salvage value'' cost structure:
\begin{eqnarray}\label{eq:smyopicii}
{\bar s}_n = \left\{ \begin{array}{ll}
             c_{n+1}-h_n, &\mbox{}n<N; \\
             s, &\mbox{}n=N.
       \end{array} \right.
\end{eqnarray}

One can interpret the modified salvage
values ${\bar s}_n $ of Eq. (\ref{eq:smyopicii}) as representing a fictitious income from  {\it inventory liquidation} (or pre-salvage at full current cost) at the beginning of the next period $n+1$, i.e., it corresponds to the situation  that the \NV\ can salvage    inventory at the price $c_{n+1}$ at the beginning of the period $n+1$. Note that the condition $c_n(1+\ell_n)+h_n\geq c_{n+1}$ is required if inventory liquidation is allowed. Otherwise, the \NV\ will stock up at an infinite level and sell them off at the beginning of period $n+1$. Such speculation is precluded  by the aforementioned condition.
Let further,
\begin{eqnarray}
{\bar a}_n&=&\frac{p_n-c_n [1+\ell_n]}{p_n-{\bar s}_n} \label{Eq:tilde-a},\\
{\bar b}_n&=&\frac{p_n-c_n [1+i_n]}{p_n-{\bar s}_n}. \label{Eq:tilde-b}
\end{eqnarray}
and the corresponding critical values are  given by
\begin{eqnarray}
{\bar \alpha}_n&=&F^{-1}_{n}({\bar a}_n) \label{Eq:myoptical-alpha}, \\
{\bar \beta}_n&=&F^{-1}_{n}({\bar b}_n). \label{Eq:myoptical-beta}
\end{eqnarray}
\noindent  For $n\in \bN,$ the   order quantity below defines   the {\it myopic
 policy} (II):
\begin{eqnarray}\nonumber
{\bar q}_n(x_n,y_n)= \left\{ \begin{array}{ll}
            ({\bar \beta}_n-x_n)^+, &\mbox{  $x_n+y_n\geq {\bar \beta}_n$}; \\[0.3cm]
             (y_n)^+, &\mbox{  ${\bar \alpha}_n\leq x_n+ y_n < {\bar \beta}_n$}; \\[0.3cm]
             {\bar \alpha}_n-x_n, &\mbox{  $x_n+y_n < {\bar \alpha}_n$.}
       \end{array} \right.
\end{eqnarray}
Let $V^{L}_{n}(x_n,y_n)$ denote the optimal expected future value when the  inventory liquidation
option is available only at the beginning of period $n+1$ (but not in the remaining periods $n+2,\ldots,N$) given the initial state  $(x_n, y_n)$ of  period $n$.
 For notational simplicity, let $\xi_{n+1}=\xi_{n+1}(x_{n},y_{n},z_{n},D_{n}) $ represent the total capital and inventory
 asset value in period $n+1$
when the \NV\ orders $z_n\geq x_n$ in state  $(x_n, y_n)$ and the demand is $D_{n}$. Note   that
$\xi_{n+1}=x_{n+1}+y_{n+1},$
where $x_{n+1}$ and $y_{n+1}$ are given by Eqs. (\ref{Eq:xn}) and (\ref{Eq:yn}),  respectively.
Therefore, $ V^{L}_n$   can be   written as
\begin{eqnarray}
V^{L}_n(x_{n},y_{n})&=&\max_{z_n\geq x_n}{\bE}\left[ V_{n+1}(0,\xi_{n+1})\,|\,x_n,y_n
\right]. \label{Eq:tildeV-e}
\end{eqnarray}

Note that for any real functions $g(x)$ and $f(x)$, with
$g(x)$ increasing and concave in $x$ and $f(x)$ concave in $x$,  the function $g\big(f(x)\big)$ is concave in $x \in \Re .$  Thus, the concavity of ${\bE}\left[ V_{n+1}(0,\xi_{n+1})\,|\,x_n,y_n\right]$  in $z_n$  follows
from the observations that $V_n(0,\xi)$ is increasing and concave  in $\xi$ (cf.  part 2 of Lemma \ref{Lem:increasing}) and
$\xi_{n+1}$ is concave in  in $z_n$, since both $x_{n+1}$ and $y_{n+1}$ are concave in $z_n$ (cf.
discussion in the proof of Lemma \ref{Lem:increasing} regarding the second order derivatives of
$x_{n+1}$ and $y_{n+1}$ with respect to $z_n$).
Therefore, $V^{L}_{n}$ has an  optimal policy characterized by a sequence of threshold value pairs, $(\alpha^{L}_{n},\beta^{L}_{n})$.

\noindent Before exhibiting  the upper bounds of $\alpha_n$ and $\beta_n$, we present the following result.
\begin{pro}\label{Pro:x-y}
For any period $n$ with initial state $(x_n, y_n)$, \\
\noindent (i) $V_n(x_n-d,\, y_n+d)$ is increasing in $d$ where $0\leq d\leq x_n$.\\
\noindent (ii) The partial derivatives satisfy, $\partial V_n(x_n,y_n)/\partial y_n\geq \partial V_n(x_n,y_n)/\partial x_n$.
\end{pro}

Part (i) of Proposition \ref{Pro:x-y} states that, for the same total asset, more capital  is more profitable than more inventory. This is true because more capital provides more flexibility of converting cash into product. Part (ii) implies that the contribution margin of a unit of cash is higher than that of a unit of inventory. In  view of Proposition \ref{Pro:x-y}, we have the following results.

\begin{pro}\label{Pro:aL}  For any state $(x_n,y_n)$,  and all $n=1,\ldots N$, the following are true:\\
(i) $V^{L}_n(x_{n},y_{n})\geq V_n(x_{n},y_{n})$,\\
(ii) $\alpha^{L}_{n} \geq \alpha_{n}$ and
  $\beta^{L}_{n} \geq \beta_{n}$ .
\end{pro}

Note that inventory liquidation at the beginning of period $n+1$ provides the \NV\ with additional flexibility, since the initial inventory $x_{n+1}$ can be liquidated into cash, thereby allowing the \NV\ to holds only $\xi_{n+1}=x_{n+1}+y_{n+1}$ in cash. Note further that the \NV\ will chose to stock up to a higher level of inventory when liquidation is allowed. Indeed,  if the \NV\ orders more in  period $n$, all leftover inventory after satisfying the period demand $D_{n}$ can be   salvaged at full cost $c_{n+1}$ at the beginning of the next period, $n+1$.  In other words, the \NV\ will take the advantage of inventory liquidation to stock product at a higher level than that corresponding to the case in which  liquidation is not allowed in the current period $n$. The advantage of doing so is twofold: (1) more demand can be satisfied leading to more revenue, and (2) no extra cost is accrued from liquidation of leftover inventory.
%


The next result establishes the upper bound properties of the myopic policy (II).
\begin{pro}\label{Pro:Optimal-myopical-tilde}
For any period $n=1,2,\ldots,N-1$,  if $c_n(1+\ell_n)+h_n\geq c_{n+1}$,  then the threshold constants of the
optimal policy given in Eqs. (\ref{Eq: N_derv profit_qi})-(\ref{Eq:
N_derv profit_ql}) and its myopic optimal policy given in Eqs.
(\ref{Eq:myoptical-alpha})-(\ref{Eq:myoptical-beta}) satisfy ${\bar \alpha}_n  \geq \alpha_n$ 
and ${\bar \beta}_n \geq  \beta_n. $\\
%
For the last period $N$, $\alpha_N={\bar \alpha}_N$ and $\beta_N = {\bar \beta}_N$.
\end{pro}

\subsection{An Algorithm to Compute $(\alpha_n,\beta_n$)}
\label{Sec:Algorithm}

It is generally not easy to compute solutions even for the simpler multi-period problems which involve base-stock levels.  In the supply chain literature, computational methods are often based on heuristic algorithms or iterative search procedures. Determination of base-stock levels in each period or for each echelon and the minimization of the total supply chain cost are notoriously complex and computationally tedious, especially for our underlying two-dimensional state space with the attending curse of dimensionality; cf.  
 \cite{Lee1993Decentral},
 \cite{Chen1994LBounds}, 
 \cite{Minner1997Algo}, 
 \cite{Shang2003SerialSC}  and
 \cite{Daniel2006heuristic}.
 
Nevertheless, with the aid of the lower and the upper bounds developed in \S \ref{Sec: Lower Bounds} and \S \ref{Sec: Upper Bounds}, respectively, we develop an efficient    algorithm for the problem herein, using  bisection. In view of Theorem \ref{Them:Dynamic_order_policy}, for each period $n\in \bN$ we aim to find the roots of $\phi(\alpha_n)=0$ and $\psi(\beta_n)=0$, where  
\begin{eqnarray}\label{Eq:phi}
\phi(\alpha_n):=\bE \left[\left( \frac{\partial V_{n+1}}{\partial x_{n+1}}-(
p'_n+h'_n)\frac{\partial V_{n+1}}{\partial y_{n+1}}\right) {\bf
1}_{\{\alpha_n>D_n\}}\right]-\left[c'_n(1+\ell_n)-p'_n\right]{\bE}\left[\frac{\partial V_{n+1}}{\partial y_{n+1}}\right] ,
\end{eqnarray}
is the first-order derivative of $G_n(z_n, x_n, y_n)$ w.r.t. $z_n$ at $z_n=\alpha_n$ for a period in which the firm took a loan; and
\begin{eqnarray}\label{Eq:psi}
\psi_n(\beta_n):=\bE \left[\left( \frac{\partial V_{n+1}}{\partial x_{n+1}}-(
p'_n+h'_n)\frac{\partial V_{n+1}}{\partial y_{n+1}}\right) {\bf
1}_{\{\beta_n>D_n\}}\right]-\left[c'_n(1+i_n)-p'_n\right]{\bE}\left[\frac{\partial V_{n+1}}{\partial y_{n+1}}\right] .
\end{eqnarray}
 is the first-order derivative of $G_n(z_n, x_n, y_n)$ w.r.t. $z_n$ at $z_n=\beta_n$  for a period in which the firm made a deposit.

\bigskip

\noindent {\bf Threshold Computation Algorithm:}
\begin{center}
\linethickness{0.3mm}
\line(1,0){450}
\end{center}
\begin{description} \small \bf
\item[STEP 0]
    \begin{tabbing}
    \qquad\qquad\quad\quad\=\kill
    \>Determine and discretize the state space of elements $x\in \mathbb{X}$ and $y\in \mathbb{Y}$; \\
    \>Specify the tolerance $\epsilon>0$;
    \end{tabbing}
\item[STEP 1]
    \begin{tabbing}
    \qquad\qquad\quad\quad\=\kill
    \>For the last period $n=N$ and  $\forall (x,y)\in \mathbb{X}\times \mathbb{Y}$, compute $V_N(x,y)$ using Eq. (\ref{Eq:V_star});
    \end{tabbing}
\item[STEP 2]
    \begin{tabbing}
    \qquad\qquad\quad\quad\=\kill
    \>Set the period index $n\leftarrow n-1$;\\
    \> Compute $\underline \alpha_n$ and $\underline \beta_n$ using Eqs. (\ref{Eq:myoptical-alpha-hat})-(\ref{Eq:myoptical-beta-hat});\\
    \> Compute $\bar \alpha_n$ and $\bar \beta_n$ using Eqs. (\ref{Eq:myoptical-alpha})-(\ref{Eq:myoptical-beta});\\
    \end{tabbing}
\item[STEP 3(a)]
    \begin{tabbing}
    \qquad\qquad\quad\quad\=\kill
    \>To compute $\alpha_n$: Set $a\leftarrow\underline\alpha_n$, $b\leftarrow \bar\alpha_n$;\\
   \> While $b-a>\epsilon$ \\
   \>   \indent Set $c\leftarrow (a+b)/2$ \\
   \>  \indent If $\phi_n(c)=0$ or $(b-a)/2<\epsilon$\\
   \>     \indent\indent Display ``solution found"; Output($\alpha_n=c$); Stop\\
    \>\indent End If\\
   \> \indent If $\phi_n(c)>0$, then $a\leftarrow c$; else $b\leftarrow c$.\\ 
   \>End While; go to STEP 3(b) to compute $\beta_n$.
    \end{tabbing}
\item[STEP 3(b)]
    \begin{tabbing}
    \qquad\qquad\quad\quad\=\kill
    \>To compute $\beta_n$: Set $a\leftarrow \underline\beta_n$, $b\leftarrow \bar\beta_n$;\\
   \> While $b-a>\epsilon$ \\
   \>   \indent Set $c\leftarrow (a+b)/2$ \\
   \>  \indent If $\psi_n(c)=0$ or $(b-a)/2<\epsilon$\\
   \>     \indent\indent Display ``solution found"; Output($\beta_n=c$); Stop\\
    \>\indent End If\\
   \> \indent If $\psi_n(c)>0$, then $a\leftarrow c$; else $b\leftarrow c$.\\ 
   \>End While; go to STEP 4.\\
    \end{tabbing}
\item[STEP 4]
    \begin{tabbing}
    \qquad\qquad\quad\quad\=\kill
    \>Compute the optimal order quantity using Eq. (\ref{Eq:q*}) for each state; \\
    \>Compute $V_n(x_n,y_{n})$ by substituting $q^*_n(x_n,y_n)$ into ${\bE}\left[V_{n+1}(x_{n+1},y_{n+1})\right]$;\\
    \>If $n>1$, go back to STEP 2;  Otherwise END.
    \end{tabbing}
\end{description}
\begin{center}
\linethickness{0.3mm}
\line(1,0){450}
\end{center}

We now discuss briefly  the complexity of the bisection search algorithm above.
This method is guaranteed to converge to the roots of $\phi_n(\alpha_n)$ and $\psi_n(\beta_n)$, since $\phi_n(\underline \alpha_n)\geq 0$ while $\phi_n(\bar \alpha_n)\leq 0$; and $\psi_n(\underline \beta_n)\geq 0$ while $\psi_n(\bar \beta_n)\leq 0$.
Specifically, if 
$c_k$ is the midpoint of the interval in the $k$-th iteration of STEP 3, then the difference between $c_k$ and a solution $c$, i.e., $\alpha_{n}$ or $\beta_{n}$, is bounded by $\displaystyle |c_k-c|\leq \frac{b-a}{2^n}$.
Therefore,
the number of iterations $k$ to attain a predetermined error (or tolerance) $\epsilon$ is given by,
$k=\log_2\left(\frac{\epsilon_0}{\epsilon}\right)=\frac{\log\epsilon_0-\log\epsilon}{\log2}$,
where $\epsilon_0$ is the range of bounds; \textit{viz.} $\epsilon_0=\bar \alpha_n-\underline \alpha_n$ for computing $\alpha_n$, and $\epsilon_0=\bar \beta_n-\underline \beta_n$ for computing $\beta_n$.

\subsection {Selling Back to the Supplier: Upper Bound for the Value Function $V_{n}(x,y)$}

The myopic policies above provide lower  and upper  bounds for $\alpha_n$ and $\beta_n$, respectively, and both of them 
provide   lower bounds for the value function $V_{n}(x,y)$. In this section we obtain upper bounds for the value 
function by considering a related problem in which the firm is given the option to sell unneeded   inventory 
back to the supplier at full price. 
Thus, for some cases, the firm has more flexibility at the beginning of a period $n\in \bN$ to readjust the inventory level as compared to just ordering from the supplier at cost $c_n$. For example, the firm is allowed to sell unneeded inventory back to the supplier or at a free-trading spot market at price $c_n$. We assume $c_{n+1}\leq c_n+h_n$ to preclude situations that the firm can profit from selling product in inventory back to the supplier. \\

In practice this situation    occurs  either when there is an    agreement between the supplier and the firm such as a \textit{vendor-managed inventory} (VMI) contract, or  there is a commodity spot market, e.g., energy sources and natural resources; cf. \cite{Secomandi2010trade}. We will  refer to such a setting as a \textit{selling back} option. The selling back  provides  the firm with the flexibility to liquidate any excess inventory into cash prior to the regular selling season. This technically reduces the state of the system from two variables, $x_{n}$ and $y_{n}$, to a single variable $\xi_n$ that represents the sum of the cash value of the on-hand inventory plus the cash value, expressed in product units, i.e., $\xi_n=x_n+y_n$. 

Adding Eqs. (\ref{Eq:xn}) and (\ref{Eq:yn}) yields the 
 dynamics equations
for $\xi_{n+1}=x_{n+1}+y_{n+1},$
with the decision variables   $z_n,$ 
\begin{eqnarray}
\xi_{n+1}=p'_{n}\cdot z_{n}-(p'_n+h'_n-1)
\left[z_{n}-D_{n}\right]^{+}
+c'_{n} \cdot
(\xi_{n}-z_{n})\bigg[(1+i_n){\bf 1}_{\{z_{n}\leq \xi_n\}}+(1+\ell_n){\bf
1}_{\{z_{n}> \xi_n\}}\bigg].
\label{Eq:xin}
\end{eqnarray}
The dynamic programming equations for the  value 
function $V^S_n(\xi_{n})$ are
\begin{equation}
V^S_n(\xi_{n})=\max_{z_n\geq 0}\bE[V^S_{n+1}(\xi_{n+1})],\label{Eq:VS} \\
\end{equation}
where the terminal value   $V^S_N(\xi_N)=V_N(0,\xi_N)$ is given by Eq. (\ref{Eq:DP-optimal-N}). \\

It is straightforward to show that $V^S_n(\xi)$ is increasing and concave in $\xi$ for any $n\in\bN$. Consequently, we have the following result.
 
\begin{thm}\label{Them:SellingBack} (\textbf{Optimal inventory-trading
policy with selling back}).\\
At the beginning of period $n\in \bN$, if the firm has a selling back option, then \\
(i) The optimal policy is ($\alpha^S_n$, $\beta^S_n$) policy, where $0\leq \alpha^S_n<\beta^S_n$. That is,
the optimal target inventory is
\begin{eqnarray}\nonumber
z^{*}_n(\xi_n)= \left\{ \begin{array}{ll}
            { \beta}^S_n, &\mbox{  $\xi_n\geq {\beta}^S_n$}; \\[0.3cm]
             \xi_n, &\mbox{  ${ \alpha}^S_n\leq \xi_n< {\beta}^S_n$}; \\[0.3cm]
             {\alpha}^S_n, &\mbox{  $\xi_n < {\alpha}^S_n$.}
       \end{array} \right.
\end{eqnarray}
(ii) For any $x_n$ and $y_n$, the value  functions satisfy
$$V^S_n(x_n+y_n)\geq V^L_n(x_n, y_n)\geq V_n(x_n,y_n).$$
\end{thm}

\bigskip

Theorem \ref{Them:SellingBack} states that, with the selling back option, the firm can freely reposition its inventory via selling down or ordering up as follows. 
If  $\xi_n\geq \beta^S_n$, then it is optimal to 
choose the  target inventory $z_n$ equal to $\beta^S_n$; 
 if    $\xi_n\leq \alpha^S_n$, then it is optimal to 
choose the  target inventory $z_n$ as $\alpha^S_n$;
 and otherwise, if $\alpha^S_n < \xi_n < \beta^S_n$,  
then it is optimal to 
choose the  target inventory $z_n$ equal to $\xi_n$.

Note that  $V^S_n$  provides an upper bound for $V_n$ which is easier 
to compute than that of $V_n(x,y),$ since the former is a function of only one
variable.

\section{Numerical Studies}
\label{Sec:numerical study}
In this section, we present numerical studies to glean further insights into the interplay between cash flow and inventory flow. To this end, we fix the following parameters across periods:
the selling price $p=2000$; 
the ordering cost $c=1000$;  
the holding cost $h=500$;   
the deposit interest rate $r=1\%$;  
and the loan interest rate $\ell=15\%$. 
The salvage value $s=600$. 
In the sequel, Subsection \ref{Subsec:NumericalStudy-N} considers the impacts of time horizon and initial states on the optimal solution; Subsection \ref{Subsec:NumericalStudy-D} studies the corresponding impact of various demand distributions, and Subsection \ref{Subsec:NumericalStudy-boundsV} investigates the performance of the lower and upper bounds of the value function. 

\subsection{Impact of Time Horizon} \label{Subsec:NumericalStudy-N}
We assume that the demand is uniformly distributed as $D\sim \mathcal{U}[0,20]$ and select $N=3$ and $N=6$ as time horizions. The solution $V_1(x,y)$, is generally obtained by backward induction using the numerical algorithm as presented in \S \ref{Sec:Algorithm},
and to emphasize the dependence of $V_1(x,y)$ on $N$ we will write $V_{1:N} (x,y)$.
Figure \ref{Fig:3d_plot_VN} depicts the resultant optimal functions $V_{1:N}(x,y)$ for the two selected time horizons and states $(x,y)$ with $x$ ranging from 0 to 20 and $y$ ranging from -50 to 50. \\

Observe first that $V_{1:N}(x,y)$ in Figure \ref{Fig:3d_plot_VN} is increasing and jointly concave in $x$ and $y$. This observation comports with Lemma \ref{Lem:increasing}.

Further, the top-right of Figure \ref{Fig:3d_plot_VN} shows that $V_{1:N}(x,y)$ is increasing in $N$ for large values of $x$ and/or $y$; however, the far bottom-left of that figure shows that for sufficiently small $x$ and/or $y$, $V_{1:N}(x,y)$ is decreasing in $N$. 
For cases of sufficiently large or small $y$, the optimal value is primarily determined by financial decisions as apposed to operational ones, e.g., there is a large cost due to loans for sufficiently negative $y$,  and consequently, $V_{1:6}(x,y)$ is smaller than  $V_{1:3}(x,y) $.  The opposite is true for sufficiently large $y$. Conversely, the middle segment of that figure shows that for non-extreme values in the range of $x$ and $y$, the operational decisions play a more important role than the financial ones, yielding significantly larger magnitudes of $V_{1:N}(x,y)$ over the time horizon  $N$. The dominant role of operational decisions is supported by the observation that the rate of increase of $V_{1:N}(x,y)$ in $y$ is increasing in $N$.

\subsection{Impact of Demand Distributions}\label{Subsec:NumericalStudy-D}
\begin{figure}[htbp]\centering
\begin{minipage}[t]{.47\textwidth}
\caption{Function $V_1(x,y)$ for $N=6$ and $N=3$ with $D\sim\mathcal{U}[0,20]$} \label{Fig:3d_plot_VN}
\includegraphics[width=0.9\textwidth]{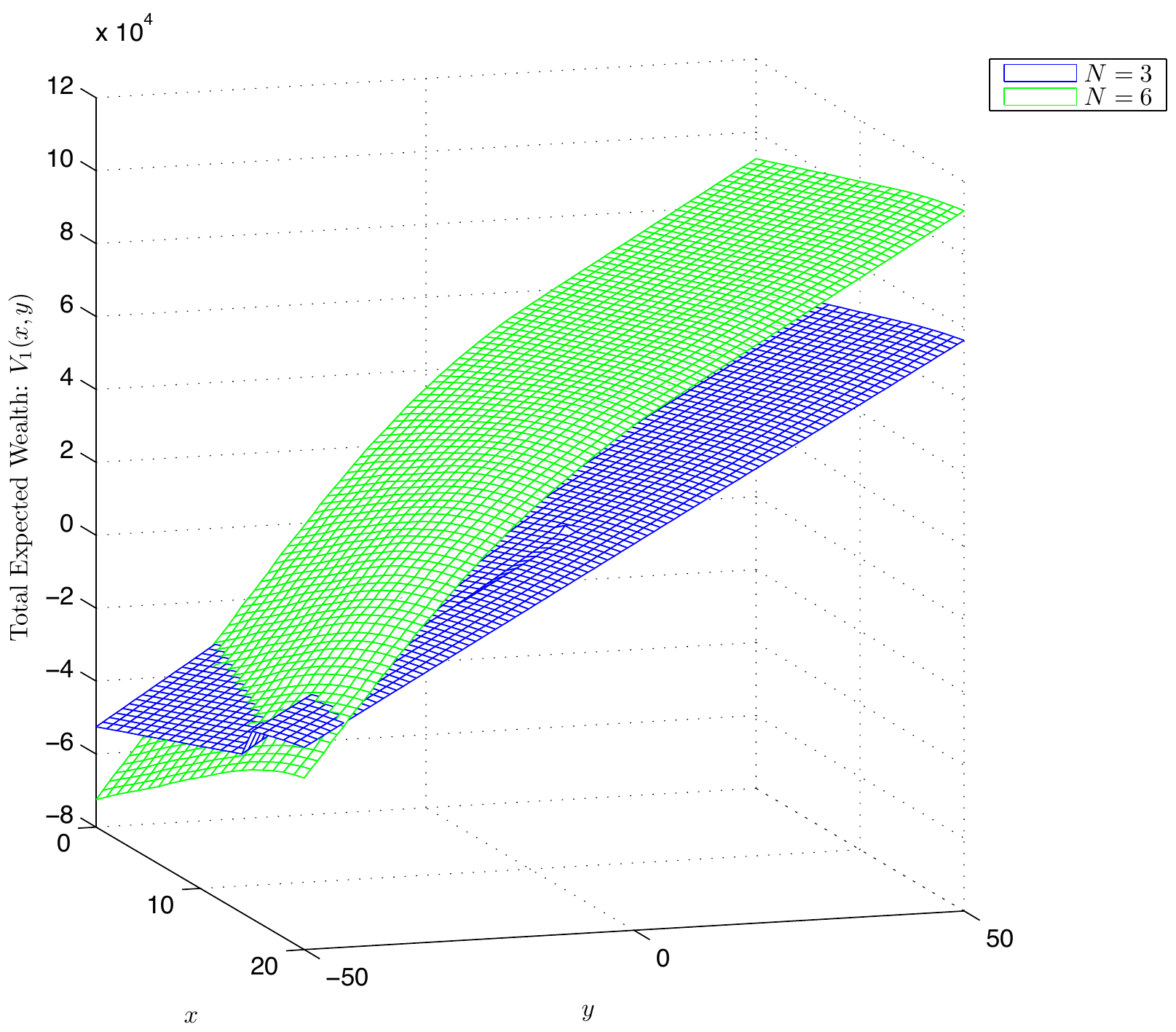}
\end{minipage}
\hfill
\begin{minipage}[t]{.47\textwidth}
\caption{Function $V_1(x,0)$ for $N=6$ with Selected Uniform Demand Distributions} \label{Fig:V_UinformD_X}
\includegraphics[width=0.9\textwidth]{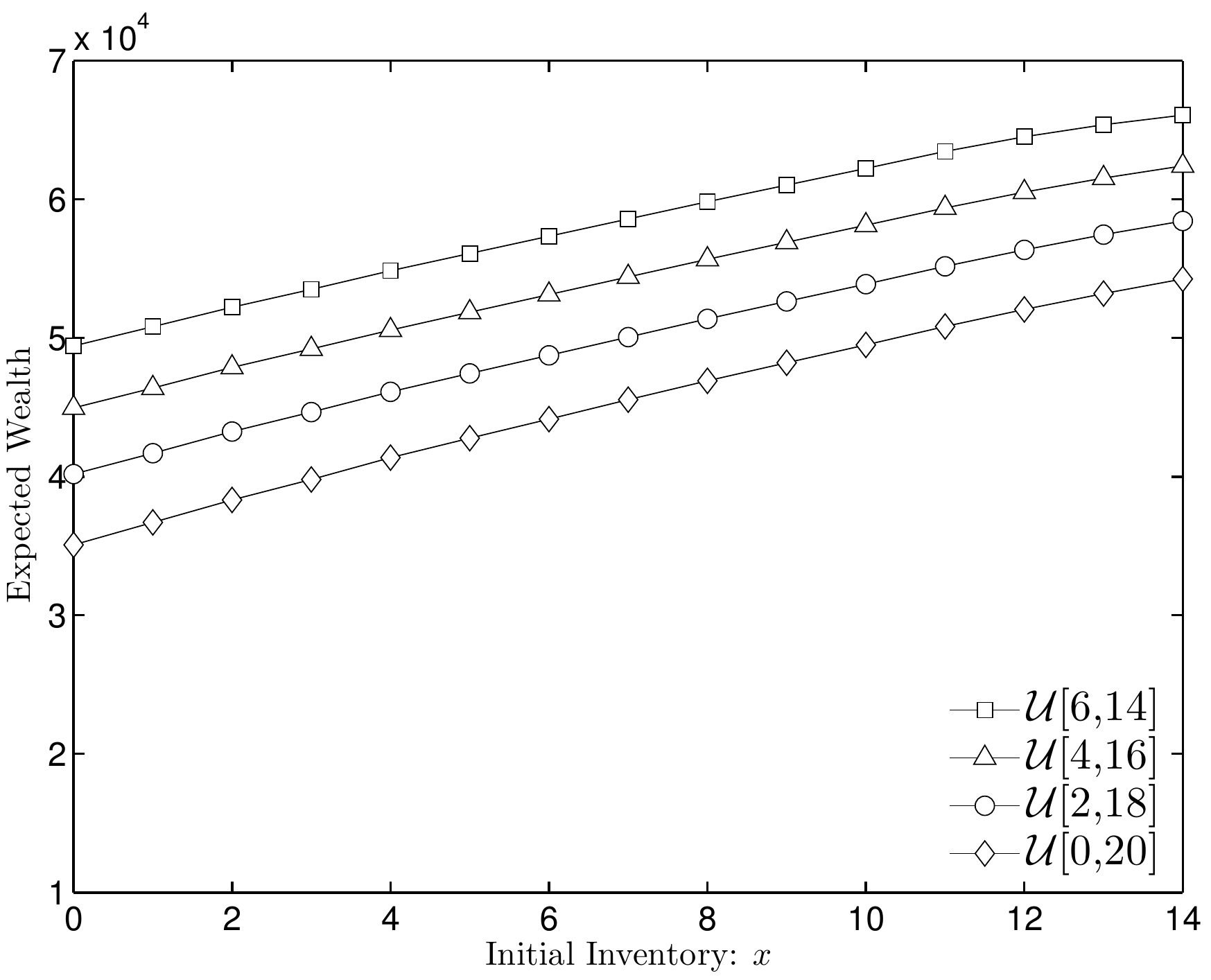}
\end{minipage}
\end{figure}

\begin{figure}[htbp]\centering
\begin{minipage}[t]{.47\textwidth}
\caption{Function $V_1(0,y)$ for $N=6$ with Selected Uniform Demand Distributions} \label{Fig:V_UinformD_Y}
\includegraphics[width=0.9\textwidth]{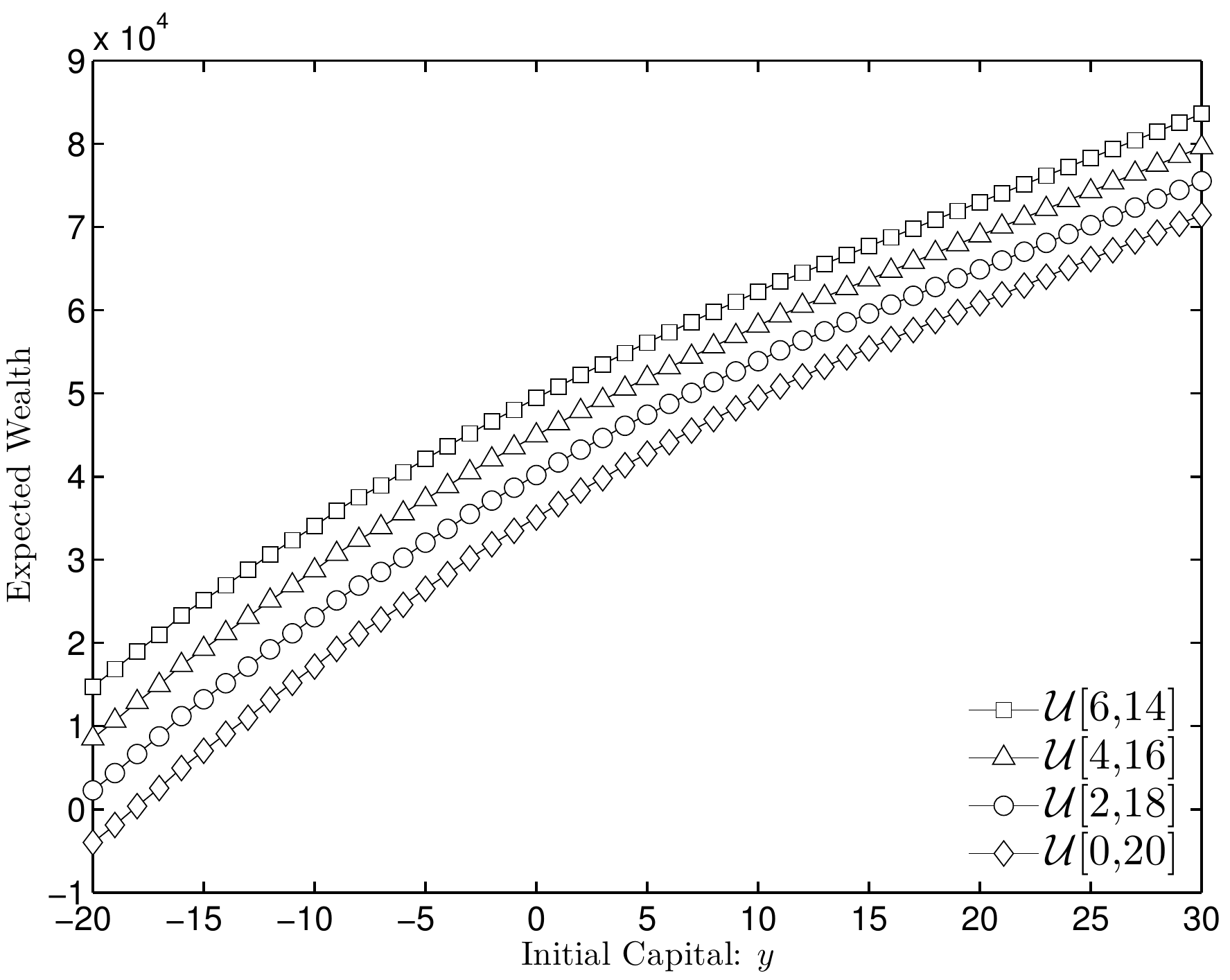}
\end{minipage}
\hfill
\begin{minipage}[t]{.49\textwidth}
\caption{$V_1(x,y)$ for $N=6$: Uniform v.s. ZIP Demand Distributions} \label{Fig:V_Uinforxm_ZIP}
\includegraphics[width=0.9\textwidth]{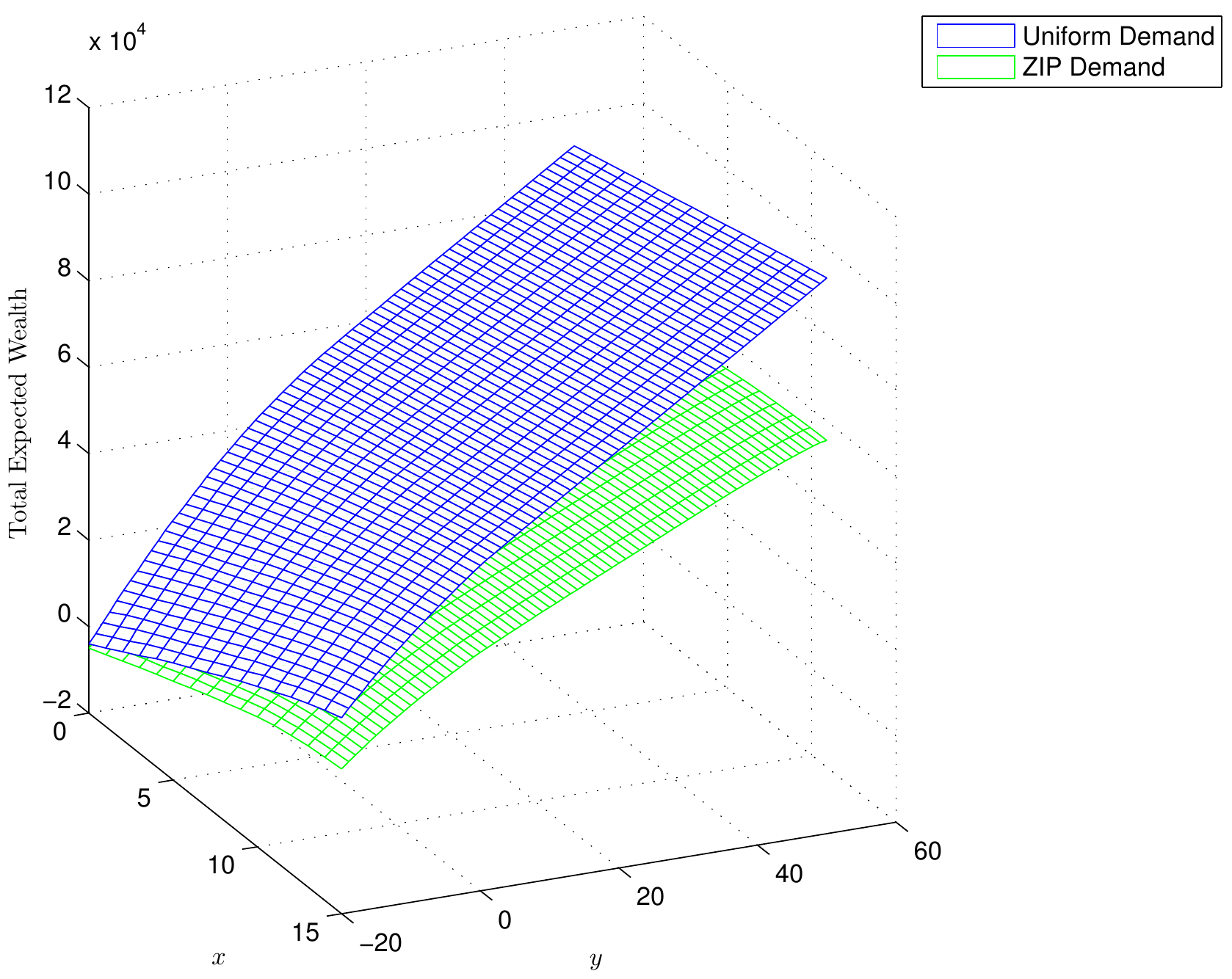}
\end{minipage}
\end{figure}

In the sequel, we fixed the time horizon as $N=6$. To study the impact of demand distributions, we first consider various uniform demand distributions with the same mean but different variance.
Figures \ref{Fig:V_UinformD_X} and \ref{Fig:V_UinformD_Y} depict $V_1(x,0)$ and $V_1(0,y)$ respectively for $D\sim \mathcal{U}[0,20]$ , $D\sim \mathcal{U}[2,18]$, $D\sim \mathcal{U}[4,16]$  and $D\sim \mathcal{U}[6,14]$. Both figures show concavity for each demand type. 
Given the common mean, the demand variance has a significant impact on the value function.  For example, there is roughly 5,000 reduction from the case of $D\sim \mathcal{U}[6,14]$ to the case of $D\sim \mathcal{U}[4,16]$ for any $x$ or $y$ in each figure.

We next consider the impact of different types of demand distributions in addition to the uniform distribution.
The \textit{Zero-Inflated Poisson} ($\mathcal{ZIP}$) distribution is a useful generalization of the Poisson distribution; cf \cite{Lambert1992ZIP}.  
The ZIP distribution with parameters $\pi\in[0,1)$ and $\lambda\geq 0$ has the following probability mass function: $\mathbb{P}(X=0)=\pi+(1-\pi)e^{-\lambda}$ and $\mathbb{P}(X=k)=(1-\pi)e^{-\lambda}\lambda^k/{k!}$ for $k=1,2,3\ldots$. Hence, its mean $\mu_X=\lambda(1-\pi)$ and its standard deviation $\sigma_X=\sqrt{\lambda(1-\pi)(1+\lambda\pi)}$. Therefore, its \textit{coefficient of variation} (defined as the ratio between the standard deviation and the mean) is $c.v.=\sqrt{\frac{1+\lambda\pi}{\lambda(1-\pi)}}$.

Figure \ref{Fig:V_Uinforxm_ZIP} compares the optimal $V_1(x,y)$ with $D\sim \mathcal{ZIP}(0.18, 10)$ and $D\sim \mathcal{U}[0,20]$, both of which have the same $c.v.=0.58$. It shows that the uniform demand yields higher $V_1(x,y)$ than the $\mathcal{ZIP}$ demand for any $x$ and $y$. This can be explained by the mean of the uniform demand ($\mu=10$) being larger than that of the $\mathcal{ZIP}$ demand ($\mu=8.2$), even though they have the same $c.v.$\\

Table \ref{Table:Profit_Gap} displays the the optimal objective function values and the performances of myopic policies (I) and (II), in terms of uniform and ZIP demands with respect to different $c.v.$ levels. In the table, for $x_1=y_1=0$, $V_1:=V_1(0,0)$, $V^{(\mathscr{X})}_1$ is the expected wealth under \textit{myopic policy} ($\mathscr{X}$), $\mathscr{X}\in \{I, II\}$, and the performance gap percentage $\Delta^{(\mathscr{X})}\%:=(V_1-V_1^{(\mathscr{X})})/V_1\times 100\%$.
First, as $c.v.$ decreases, the value of $V_1(0,0)$ increases significantly for each type of demands.
Second, under the uniform demand, myopic policy (II) with $\Delta^{(II)}\%<0.2\%$ performs significantly better than myopic policy (I) with $\Delta^{(I)}\%>3.5\%$; but the converse is observed under the $\mathcal{ZIP}$ demand. 
For the same $c.v.$ value, $\mathcal{ZIP}$ demand has lower $V_1$ than Uniform demand; this can be explained by the smaller mean of $\mathcal{ZIP}$ demand, and the significantly higher probability of zero demand, e.g., when $c.v.=0.46$, $\bP(D=0)=9\%$ for $\mathcal{ZIP}$ demand, but $\bP(D=0)=0\%$ for the Uniform demand.
Finally, note that for the last row in Table \ref{Table:Profit_Gap}, $\mathcal{ZIP}(0,10)$ is the standard Poisson distribution with mean $\lambda=10$ - the same mean as $D\sim \mathcal{U}(6,14)$. However, $\mathcal{ZIP}(0,10)$ has a smaller $V_1$, because it has a larger standard deviation, and hence a higher $c.v.$ value.

\tabcolsep=0.11cm
\begin{table}[htbp]\centering \caption{The Optimal v.s. Myopic Policies: Uniform and ZIP Demand with Various $c.v.$ ($x_1=y_1=0$ and $N=6$)}\label{Table:Profit_Gap}
\scriptsize
\begin{tabular}{ll r r lr r lr   rr  lr r r lr r lr}
\hline\hline\\[0.1ex]
\multicolumn{9}{c}{$D\sim\mathcal{U}[a,b]$}&&&\multicolumn{9}{c}{$D\sim\mathcal{ZIP}(\pi,\lambda)$} \\[0.1ex]
 \cline{1-9}\cline{12-20}\\[0.1ex]
     &      &Optimal  &&\multicolumn{2}{c}{Myopic I} &&\multicolumn{2}{c}{Myopic II}&&&
     &      &Optimal  &&\multicolumn{2}{c}{Myopic I} &&\multicolumn{2}{c}{Myopic II}\\
     \cline{3-3}\cline{5-6} \cline{8-9}
     \cline{14-14}\cline{16-17} \cline{19-20}\\[0.08ex]
 $c.v.$ &  $[a,b]$  &$V_1$ &  &$V^{(I)}_1$ &$\Delta^{(I)}$\%  & &$V^{(II)}_1$ &$\Delta^{(II)}$\%  & & &
 $c.v.$ &  $(\pi,\lambda)$  &$V_1$ &  &$V^{(I)}_1$ &$\Delta^{(I)}$\%  & &$V^{(II)}_1$ &$\Delta^{(II)}$\% \\[0.08ex]
  \cline{1-9}\cline{12-20}\\[0.1ex]
0.58&[0, 20]&	35074 &&	30271&	13.69\% &&	35016&	0.16\%	&&&	0.58&	(0.18,	10)&	23130&&	22800&	1.43\%&&	21757&	5.94\%\\[0.2ex]
0.46&[2,18] &	40174 &&	36784&	8.44\%	&&  40158&	0.04\%	&&&	0.45&	(0.09,	10)&	26920&&	26910&	0.04\%&&	26142&	2.89\%\\[0.1ex]
0.35&[4,16] &	44950 &&	42329&	5.83\%	&&  44886&	0.14\%	&&&	0.35&	(0.02,	10)&	29612&&	29484&	0.43\%&&	29115&	1.68\%\\[0.1ex]
0.23&[6,14] &	49428 &&	47677&	3.54\%	&&  49385&	0.09\%	&&&	0.31&	(0,	10)    &	30355&&	30288&	0.22\%&&	29911&	1.46\%\\[0.1ex]
\hline
\end{tabular}
\end{table}

\begin{figure}[htbp]\centering
\begin{minipage}[t]{.47\textwidth}
\caption{$V^S_1(\xi)$ v.s. $V_1(0,y)$ and $V_1(x,0)$ for $N=6$ with $D\sim\mathcal{U}[0,20]$} \label{Fig:Fig_V_UinformD_comparing}
\includegraphics[width=0.9\textwidth]{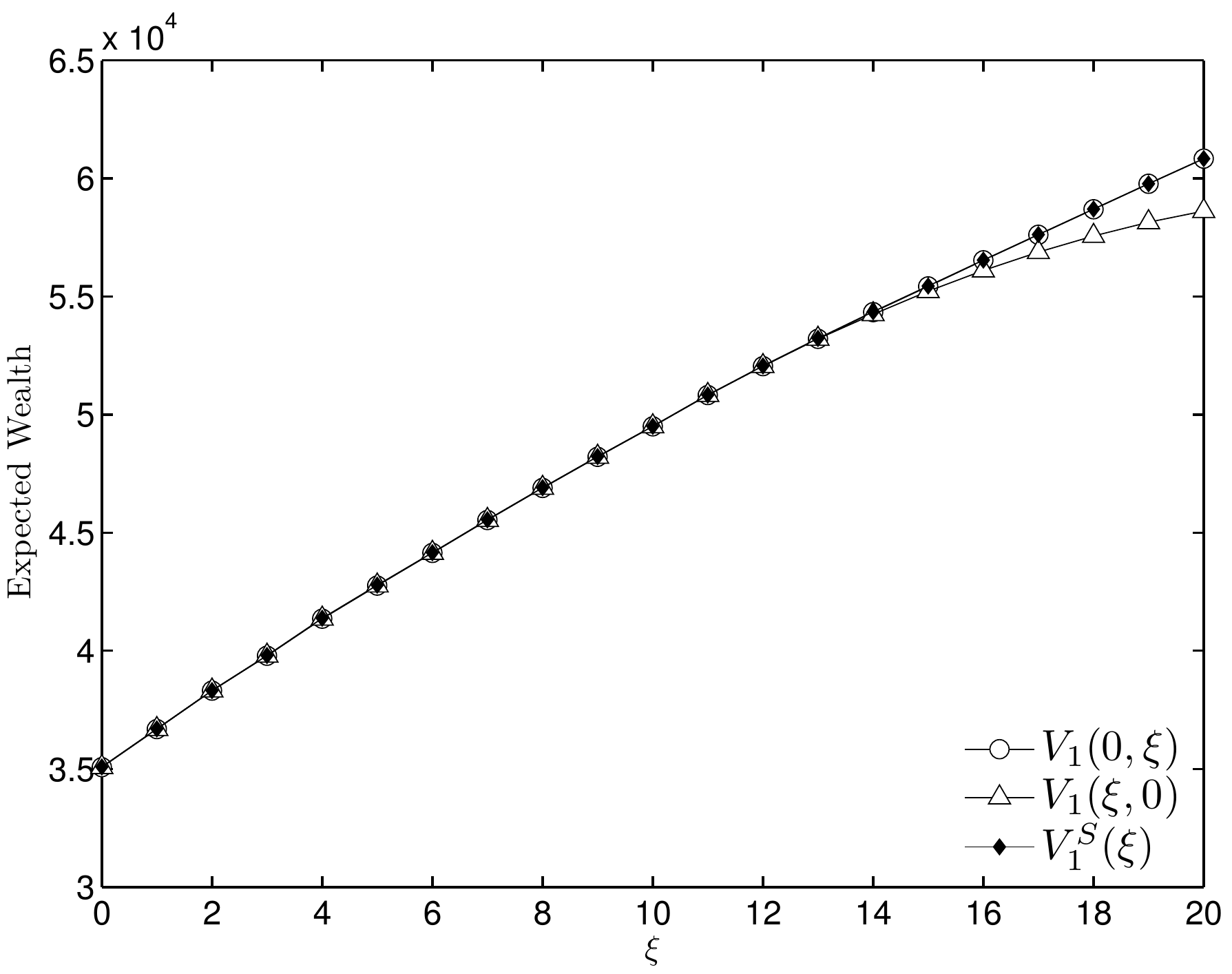}
\end{minipage}
\hfill
\begin{minipage}[t]{.47\textwidth}
\caption{Function $V^S_1(\xi)$ for $N=1,2,4,6$ with $D\sim\mathcal{U}[0,20]$} \label{Fig_V_UinformD_V_S(N)}
\includegraphics[width=0.9\textwidth]{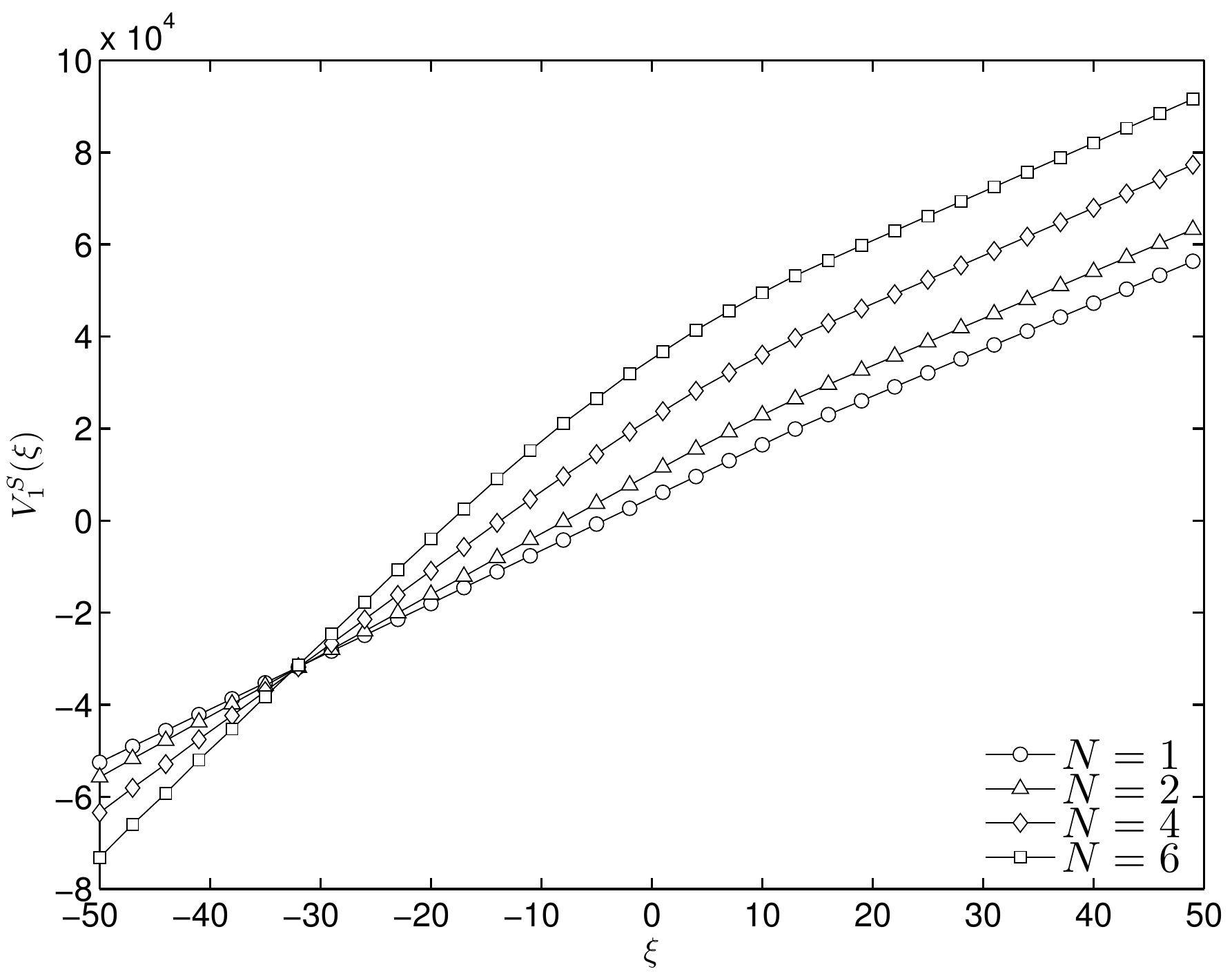}
\end{minipage}
\end{figure}

\subsection{Upper and Lower Bounds of the Value Function}\label{Subsec:NumericalStudy-boundsV}
As discussed in \S \ref{Sec:myopic}, we know that each of the myopic policies provides a lower bound for the value function, hence we can take the larger one as the better lower bound. Furthermore, $V^{{S}}_{1}$ provides an upper bound. In what follows, we numerically examine their approximating performance  with  the same parameters as the ones selected in the previous numerical studies.

Figure \ref{Fig:Fig_V_UinformD_comparing} compares $V^{S}_{1}(\xi)$  with $V_{1}(\xi,0)$ and $V_{1}(0,\xi)$ for $N=6$. First, for any $\xi$, it is shown that $V^{S}_{1}(\xi)\geq V_{1}(0,\xi)\geq V_{1}(\xi,0)$. Second, the difference between $V^{S}_{1}(\xi)$ and $ V_{1}(0,\xi)$ is almost imperceptible which implies that the upper bound is fairly tight for the state of high capital level $y$ (or low product level $x$). Further, the bound gap between $V^{S}_{1}(\xi)$ and $V_{1}(\xi,0)$ gets larger as $\xi$ increases.

Figure \ref{Fig_V_UinformD_V_S(N)}  depicts the effect of $N$ on  $V^{S}_{1:N}(\xi)$, for $N=1,2,4,6$. Observe that for any $N$, $V^{S}_{1:N}(\xi)$ is increasing and jointly concave in $\xi$.  Further, referring to the top-right of Figure \ref{Fig_V_UinformD_V_S(N)}, $V^{S}_{1:N}(\xi)$ is increasing in $N$ for large values of $\xi$, e.g., $\xi>10$; however, referring to the far bottom-left of the figure, $V^{S}_{1:N}(\xi)$ is decreasing in $N$ for sufficiently small $\xi$, e.g., $\xi<-30$. 
Finally, we conclude that operational decisions plays a more important role than the financial decision for the middle range of $\xi$, e.g., $\xi\in (-30,10)$, since the slope at each $N$ in the interval is significantly larger than elsewhere, and gets larger as $N$ increases.\\

Table \ref{Table:V_Bounds} exhibits the
\textit{optimal value function} ($\mathrm{Opt.}$) and its  lower and upper bounds for
$D\sim \mathcal{U}[0,20]$ and $D\sim \mathcal{U}[6,14]$, $N=6, 12$ and selected product level state $x=0,7,14$. In this case, the lower bound is obtained from myopic policy (II) since it performs better than myopic policy (I).
%
First, each bound provides a fairly tight approximation for the optimal value with $\Delta^{(\mathscr{X})}\%$ less than 1\%.
Second, as $c.v.$ of demand distribution decreases from $\mathcal{U}[0,20]$ to $\mathcal{U}[6,14]$, the optimal value function and its lower and upper bounds increase, due to the reduced volatility of the demand.  Further, the lower bound gap $\Delta^{(II)}$ gets smaller, while this change is not evident for its upper bound gap $\Delta^{S}$.
Third, for each $N$ and each demand distribution, as $x$ increases, the bound gap $\Delta^{S}$ gets larger, whereas this behavior is not evident  for $\Delta^{(II)}$.
Finally, as the time horizon $N$ gets longer, the optimal value function and its lower and upper bounds get larger, and the bound gap 
$\Delta^{(\mathscr{X})}$ gets larger as well for both bounds $\mathscr{X}\in \{II, S\}$. However, the change of the corresponding 
gap percentage $\Delta^{(\mathscr{X})}\%$ is not evident.

\begin{table}[htbp]\centering \caption{The Lower and Upper Bounds for the Value Function $V_{1}(x,0)$}\label{Table:V_Bounds}\scriptsize
\begin{tabular}{l ll   l rcl c rcl    cc     l rcl c rcl }
\hline\hline\\[0.1ex]
 & &\multicolumn{8}{c}{$D\sim\mathcal{U}[0,20]$}&&&\multicolumn{8}{c}{$D\sim\mathcal{U}[6,14]$} \\[0.1ex]
 \cline{3-10}\cline{13-20}\\[0.1ex]
    &   &$\mathrm{Opt.}$   & \multicolumn{3}{c}{$\mathrm{Lower Bound}$}  &     & \multicolumn{3}{c}{$\mathrm{Upper Bound}$}   &           &     & $\mathrm{Opt.}$   &\multicolumn{3}{c}{$\mathrm{Lower Bound}$}  & & \multicolumn{3}{c}{$\mathrm{Upper Bound}$}   \\
\cline{3-3}\cline{4-6}\cline{8-10}\cline{13-13}\cline{14-16}\cline{18-20}\\
    $N$         &  $x$        & $V_{1}$  &$V^{(II)}_{1}$ &$\Delta^{(II)}$&$\Delta^{(II)}$\%   & & $V^{S}_{1}$&$\Delta^{S}$& $\Delta^{S}$\%  &  & & 
                       $V_{1}$  &$V^{(II)}_{1}$&$\Delta^{(II)}$ &$\Delta^{(II)}$\%   & & $V^{S}_{1}$&$\Delta^{S}$&$\Delta^{S}$\% \\
\cline{1-2}\cline{3-3}\cline{4-6}\cline{8-10}\cline{13-13}\cline{14-16}\cline{18-20}\\[0.1ex]
\multirow{3}{*}{\textbf{6}}&   0	&	35074&	35016&	58&	0.16\%&		&35080&	-6&	-0.02\%&	&&	 49428&	49386&	43&	0.09\%&	&	49428&	0&	0.00\%\\
&  7 &	45542&	45435&	107&	0.24\%&		&45550&	-9&	-0.02\%&	&&	58575&	58536&	39&	0.07\%&	&	58575&	0&	0.00\%\\
&14 &	54248&	54200&	48&	0.09\%&		&54355&	-107&-0.20\%&	&&	66076&	66036&	40&	0.06\%&	&	66634&	-558&	-0.84\%\\[0.2ex]
\hline\hline\\

\multirow{3}{*}{\textbf{12}}
&0	&75888&		75693&	194&	 0.26\%&		&75923&	-36&	-0.05\%&	&&	103872&		103760&	112&	 0.11\%&	&	103872&	0&	0.00\%\\
&7	&87273&		87057&	216&	  0.25\%&		&87290&	-18&	-0.02\%&	&&	113564&		113464&	100&	 0.09\%&	&	113564&	0&	0.00\%\\
&14	&96528&		96410&	118&  0.12\%&		&96660&	-132& -0.14\%& &&	121521&		121419&	102& 0.08\%&	&	122114&	-592&-0.49\%\\
\hline\hline
\end{tabular}
\end{table}

\section{Model Extensions}
\label{Sec:discussion}
In this section, we study three extensions of the basic model. The first extension is a piecewise linear structure of the interest rate functions for both loans and deposits,  the second is the existence of a maximum loan limit and the third is backorder of unmet demand. 
\subsection {Piecewise Type of Loan and Deposit Functions}
So far,  the loan function was assumed to be a linear function:
$L(x)=(1+\ell)\cdot x$ with a constant  loan rate $\ell$. However, $L(x)$   can have a  more
complex  form in practice. In this section we investigate the often occurring case in which $L(x)$   is a piecewise linear function, i.e., it has the form,
\begin{eqnarray}
 L(x)=(1+\ell^{(m)}) \cdot x, \qquad x\in(x^{(m-1)},x^{(m)}],\nonumber
 \end{eqnarray}
where $x^{(m-1)}<x^{(m)}$, $x^{(0)}=0$ and $\ell^{(m)}< \ell^{(m+1)}$ for
$m=1,2,3,...$.

Similarly, we consider  the deposit interest function to be a piecewise linear function  of the form,
\begin{eqnarray}
M(y)=(1+i^{(k)}) \cdot y, \qquad y\in(y^{(k-1)},y^{(k)}],\nonumber
 \end{eqnarray}
where $y^{(k-1)}<y^{(k)}$, $y^{(0)}=0$ and $i^{(k)}\leq i^{(k+1)}$ for
$k=1,2,3,...$.

Without loss of generality, we assume that the loan interest rates are always greater  than the deposit interest rates, that is ${\bar i}< \ell^{(1)} $ where ${\bar i}=\sup_k \{i^{(k)}\}$.

\noindent Before characterizing the optimal ordering policy, we  introduce the threshold values of $\alpha^{(m)}$ and
$\beta^{(k)}$ for $m,k=1,2,3,...$ as follows
\begin{eqnarray}
\alpha^{(m)}&=&F^{-1}(a^{(m)}) ; \label{Eq:alpha_m_piece}\\
\beta^{(k)}&=&F^{-1}(b^{(k)}) , \label{Eq:beta_k_piece}
\end{eqnarray}
\noindent where
$$a^{(m)}=\frac{p-c\cdot(1+\ell^{(m)})}{p-s} ; \,\, b^{(k)}=\frac{p-c \cdot (1+i^{(k)})}{p-s}. $$

%
It is straightforward to see that
$$\beta^{(1)}\geq\cdots \beta^{(k)}\geq \beta^{(k+1)}\cdots \geq {\bar \beta} > \alpha^{(1)}\geq\cdots \alpha^{(m)}\geq \alpha^{(m+1)}\cdots \geq 0,$$ where $\displaystyle {\bar \beta}=F^-(\frac{p-c \cdot (1+{\bar i})}{p-s})$.
We present the  structure of the optimal ordering policy in Theorem \ref{Corr:Optimal-q-piecewise}
below.

\begin{thm}\label{Corr:Optimal-q-piecewise}
For any given initial inventory-capital state $(x,y)$, the optimal order quantity is
\begin{eqnarray}\nonumber
q^{*}(x,y)= \left\{ \begin{array}{ll}
           (\beta^{(k)}-x)^{+}, &\mbox{  $ \beta^{(k+1)} \leq x+y\leq \beta^{(k)} $}; \\
            \cdots, & \cdots\\
            (y)^{+}, &\mbox{  $\alpha^{(1)}\leq x+ y < {\bar \beta}$}; \\
            \cdots, & \cdots\\
  (\alpha^{(m)}-x)^{+}, &\mbox{  $\alpha^{(m+1)}\leq x+y < \alpha^{(m)}$,}
       \end{array} \right.
\end{eqnarray}
where $\{\alpha^{(m)}\}$ and $\{\beta^{(k)}\}$ are given by Eqs. (\ref{Eq:alpha_m_piece}) and
(\ref{Eq:beta_k_piece}), respectively.
\end{thm}

For the multi-period problem, the optimal order quantity for each period has a structure similar to that of its single period counerpart, and the proof is analogous to that of Theorem \ref{Corr:Optimal-q-piecewise}. We shall note that for period $n<N$, the threshold values, $\alpha_{n}^{(m)}$ and $\beta_{n}^{(k)}$, are determined by $\xi_{n}$. However, for the last period $N$, the threshold values, $\alpha_{N}^{(m)}$ and $\beta_{N}^{(k)}$, are independent of the initial state $x_{N}$ and $y_{N}$, or $\xi_{N}$.
%
 We shall mention that if the lender charges a high loan rate such that $p\leq c\cdot (1+\ell^{(m)})$ for loan amount larger than $x^{(m-1)}$, then $\alpha^{{(m)}}=0$ via Eq. (\ref{Eq:alpha_m_piece}), since $a^{(m)}$ is nonpositive. In view of Theorem \ref{Corr:Optimal-q-piecewise}, if the total asset $x+y<0$, then the firm will not order via a loan.

\subsection {Financing Subject to a Maximum Loan Limit Constraint}
In practice, the outstanding loan amount is often restricted not to exceed a maximum limit.   Let $L_n>0$ denote the maximum loan limit   for period $n$; hence the maximum loan can finance finance $L'_n:=L_n/c_n$ units of the product. In this case, we have the following structural results for the optimal ordering policy.
\begin{thm}\label{Them:Dynamic_order_policy_Loan_Capacity} (\textbf{Optimal ordering
policy under a maximum loan limit}).\\
For period $n\in \bN$ with given state $(x_n,y_n)$ at the
beginning of the period, if there is a loan limit $L_n$, then there exist positive constants
$\alpha^L_n=\alpha^L_n(\xi_n)$ and $\beta^L_n=\beta^L_n(\xi_n)$ with
$\alpha^L_n\leq \beta^L_n$, which characterize the  optimal order quantity as
follows:
\begin{eqnarray}\nonumber
q^{*}(x_n,y_n)= \left\{ \begin{array}{ll}
            (\beta^L_n-x_n)^+, &\mbox{  $x_n+y_n\geq \beta^L_n$}; \\[0.3cm]
             (y_n)^{+}, &\mbox{  $\alpha^L_n\leq x_n+ y_n < \beta^L_n$}; \\[0.3cm]
             (\alpha^L_n-x_n)^{+}, &\mbox{  $\alpha^L_n-L'_n \leq x_n+y_n < \alpha^L_n$};\\[0.3cm]
             L'_n, &\mbox{  $ x_n+y_n < \alpha^L_n-L'_n$.}
       \end{array} \right.
\end{eqnarray}
\end{thm}

Note that for period $n$, the inequality $\alpha^L_n<L'_n$ is admissible. In this case, the maximum loan limit is not binding since the optimal order quantity satisfies $q^*_n\leq \alpha^L_n <L'_n$.

\subsection {Backorder System}\label{Sec:Ext:Backorder}
Recall that our basic model assumes the lost-sale rule for unmet demand. This section considers the case in which backorders are allowed and shows that  an  ordering policy of the same $(\alpha_n, \beta_n)$ structure is  optimal. 

Suppose that any unmet demand is backlogged subject to $b_n$ penalty cost per unit of unmet demand. Let ${\tilde V}_n(x,y)$ denote the optimal expected profit function of period $n$ for given initial state $(x,y)$. Accordingly, we have the following dynamic programming formulation:
\begin{eqnarray}
{\tilde V}_n(x_{n},y_{n})= \sup_{z_n \geq x_n}{\bE}\bigg[{\tilde V}_{n+1}(x_{n+1},y_{n+1}) |x_{n},y_{n}\bigg], \hspace{20 pt}
n=1,2,\cdots,N-1 \label{Eq:DP-optimal-bo}
\end{eqnarray}
where  the expectation is taken with respect to $D_n$, and $x_{n+1} ,$  $y_{n+1}$ are respectively given by
\begin{eqnarray}
x_{n+1}&=&z_{n}-D_{n}; \label{Eq:xn0-backorder}\\
y_{n+1}&=&
[{\tilde R}_{n}(D_n, z_{n})+{\tilde K}_{n}(z_{n},\xi_{n})]/c_{n+1},
\label{Eq:yn0-backorder}
\end{eqnarray}
and
\begin{eqnarray}
{\tilde R}_{n}(D_n,z_{n})
&=&(p_{n}+b_n)\cdot z_{n}-(p_n+h_n+b_n)
\left[z_{n}-D_{n}\right]^{+}-b_n\bE[D_n];\label{Eq:r-backorder}\\
 {\tilde K}_{n}(z_{n},\xi_{n})&=&c_{n} \cdot
(\xi_{n}-z_{n})\left[(1+i_n){\bf 1}_{\{z_{n}\leq \xi_n\}}+(1+\ell_n){\bf
1}_{\{z_{n}> \xi_n\}}\right]. \label{Eq:k-backorder}
\end{eqnarray}
For the terminal period $N,$
${\tilde V}_N(x_{N},y_{N})= \sup_{q_N \geq 0}{\bE}\bigg[
{\tilde R}_{N}(D_N,z_{N})+{\tilde K}_{N}(z_{N},\xi_{N}) \bigg],
$
where $h_N=-s$.


\begin{thm}\label{Them:Dynamic_order_policy-backorder} 
Under the backorder rule, for any period $n\in \bN$ with initial state $(x_n,y_n)$, there exist positive constants
$\alpha_n=\alpha_n(\xi_n)$ and $\beta_n=\beta_n(\xi_n)$ satisfying 
$\alpha_n\leq \beta_n$, which characterize the optimal order quantity as follows:
\begin{eqnarray}\label{Eq:q*_backorder}
{\tilde q}^{*}(x_n,y_n)= \left\{ \begin{array}{ll}
            (\beta_n-x_n)^+, &\mbox{  $x_n+y_n\geq \beta_n$}; \\
             (y_n)^+, &\mbox{  $\alpha_n\leq x_n+ y_n < \beta_n$}; \\
             (\alpha_n-x_n)^{+}, &\mbox{  $x_n+y_n < \alpha_n$.}
       \end{array} \right.
\end{eqnarray}
Further,  ${\alpha}_n$ is uniquely identified by
\begin{eqnarray}\label{Eq: N_derv profit_qi_backorder}
\bE \left[\left( \frac{\partial {\tilde V}_{n+1}}{\partial x_{n+1}}-(
p'_n+h'_n)\frac{\partial {\tilde V}_{n+1}}{\partial y_{n+1}}\right) {\bf
1}_{\{\alpha_n>D_n\}}\right]=\left[c'_n(1+\ell_n)-p'_n\right]{\bE}\left[\frac{\partial {\tilde V}_{n+1}}{\partial y_{n+1}}\right] ,
\end{eqnarray}
and ${\beta}_n$ is uniquely identified by
\begin{eqnarray}\label{Eq: N_derv profit_ql_backorder}
\bE \left[\left( \frac{\partial {\tilde V}_{n+1}}{\partial x_{n+1}}-(
p'_n+h'_n)\frac{\partial {\tilde V}_{n+1}}{\partial y_{n+1}}\right) {\bf
1}_{\{\beta_n>D_n\}}\right]=\left[c'_n(1+i_n)-p'_n\right]{\bE}\left[\frac{\partial {\tilde V}_{n+1}}{\partial y_{n+1}}\right] .
\end{eqnarray}
where the expectations are taken with respect to $D_n$ conditionally on the initial state $(x_n,y_n)$.
\end{thm}

\section{Concluding Remarks}
\label{Sec:conclusions}

In this paper, we studied the optimal operational-financial policy for a single-item inventory system that adimitted  both interest-bearing loans and interest-earning deposits. We showed that the optimal ordering policy, called $(\alpha_n,\beta_n)$-policy, for each period is
characterized by two threshold variables.
we presented two myopic policies which provided each a lower and upper bound for the threshold variables. The two bounds were then employed to develop an algorithm to compute the two threshold values, $\alpha_n$ and $\beta_n$.
The two myopic policies provided lower bounds of the value function. However, we further introduced an upper bound for the value function via assuming that the firm can sell unneeded inventory back to the supplier.

The results of this paper directly suggest several directions for further research.  These include  
i) addition of a fixed ordering cost and/or a fixed financial transaction cost;
ii)  addition of a fixed system operating cost proportional to the elapsed time; and
iii) incorporation of bankruptcy risk. Our study assumes risk neutral decision making. 
It would be of interest to analyze the risk-induced performance of the system, e.g., credit rollover risk and bankruptcy probabilities; cf. \cite{babich2012risk}, \cite{Gong2014dynamic} and \cite{Li2013control}.



 
 \bigskip

\noindent \textbf{Acknowledgement}

\noindent 
The authors are grateful to Dr. John Birge (the Departmental Editor), the Associate Editor and the reviewers for valuable comments which significantly improved the quality and presentation of the paper. The authors are also indebted to the valuable comments and useful suggestions provided by
Dr. Matthew Sobel 
and Dr. Xiuli Chao.
 


\bibliographystyle{ormsv080}
\bibliography{nvref}

\newpage
\section{Appendix}\label{Sec:Appendix}
\textbf{Proof of Lemma \ref{Lem:concave}.}
The continuity follows immediately from Eq. (\ref{Eq: E
profit}). We next prove the concavity via examining its first-order
and second-order derivatives.
To this end, differentiating Eq.
(\ref{Eq: E profit}) via Leibniz's integral rule yields
%
%
%
\begin{eqnarray}\label{Eq: derv profit}
\frac{\partial}{\partial q}G(q,x,y) &= &\left\{ \begin{array}{rl}
 p-c\cdot (1+i)-(p-s)F(q+x) &\mbox{ if $q<y$} , \\[0.3cm]
 p-c \cdot (1+\ell)-(p-s)F(q+x) &\mbox{ if $q> y$}
  .     \end{array} \right.
\end{eqnarray}
Therefore, for $q>y$ or $q<y$
\begin{eqnarray}
\frac{\partial^2}{\partial q^2 }G(q,x,y) =-(p-s)f(q+x).
\label{Eq:Pi_2_deriv}
\end{eqnarray}
\noindent The concavity  in $q$  now readily follows since
$\frac{\partial^2}{\partial q^2}G(q,x,y) \leq 0$ by Eq. (\ref{Eq:Pi_2_deriv}).

Although the value function is not first order and/or second order differentiable at certain points (e.g., at $q=y$), we can still consider its derivatives to show its increasing or decreasing properties, which allows us to study the optimal solution. Such consideration will be used throughout the paper.

The increasing property of   $G(q,x,y)$ in $x$ and $y$ can be
shown by taking the first order derivatives
\begin{eqnarray}
\frac{\partial}{\partial x}G(q,x,y)
&= &p{\bar F}(q+x)+sF(q+x)>0 \label{Eq: derv profit_x}  , \\
\frac{\partial}{\partial y}G(q,x,y)
&= &c \cdot \left[(1+i){\bf
1}_{\{q< y\}}+(1+\ell){\bf 1}_{\{q> y\}}\right]>0 \nonumber. 
\end{eqnarray}
The joint concavity of $G(q,x,y)$ in $x$ and $y$ can be
established  by computing  the second order derivatives below using again  Eq. (\ref{Eq:
E profit}).
\begin{eqnarray}
\frac{\partial^2}{\partial x^2}G(q,x,y) &= &-(p-s)f(q+x)<0 \nonumber ,\\
\frac{\partial^2}{\partial y^2}G(q,x,y) &= & 0 \nonumber ,\\
\frac{\partial^2}{\partial x \partial y}G(q,x,y) &= & 0 \nonumber.
\end{eqnarray}
Thus the Hessian matrix is negative semi-definite and the proof is complete.  
%
 

%
\begin{figure}[b]
  \centering
   \includegraphics[width=0.7\textwidth]{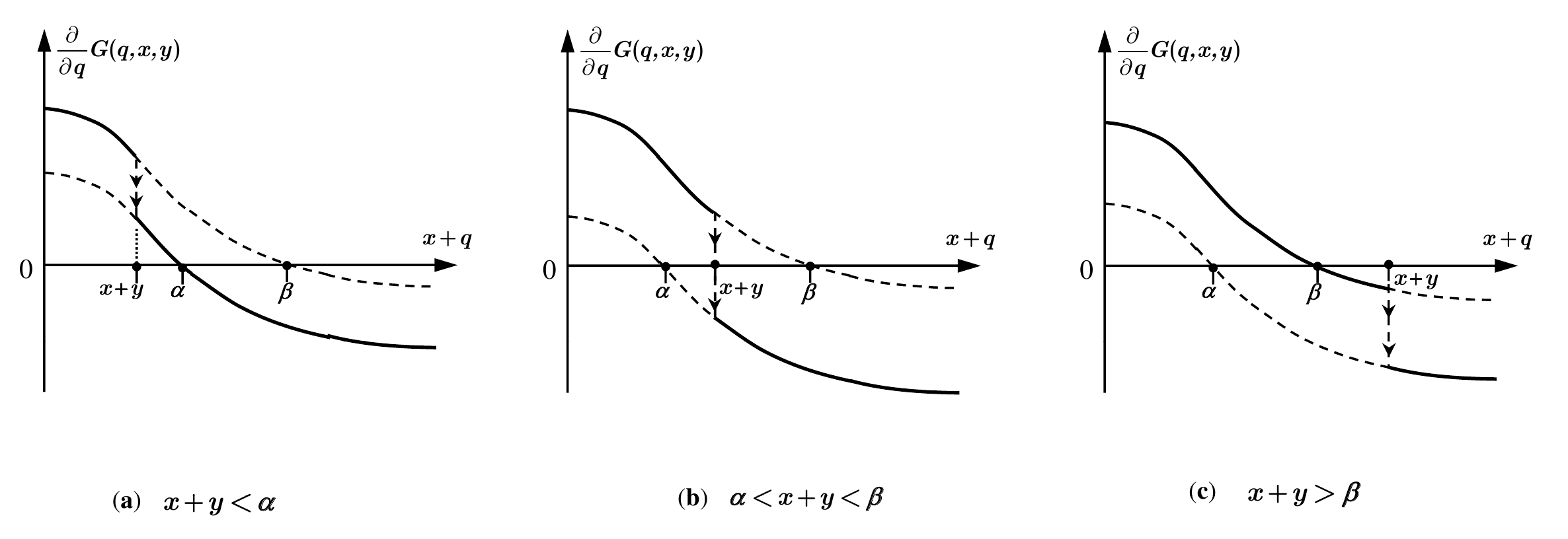}
\caption{Functional Structure for the Derivative of $G(q,x,y)$ with
Respective to $q$ } \label{Fig:derivative}
\end{figure}

\noindent\textbf{Proof of Theorem \ref{Them:Optimal-v} .}
Part (i) follows by substituting $q^*$ given by Eq. (\ref{eqstar}) into Eq. (\ref{Eq: E profit}).

For part (ii) the increasing property of $V$ can be readily justified. To show the concavity of $V$, note that by Lemma \ref{Lem:concave}, $G(q,x,y)$ is concave in
$q$, $x$ and $y$. Maximizing $G$ over $q$ and using Proposition A.3.10 of \cite{zipkin2000foundations}, p. 436, 
we have that  the
concavity in $x$ and $y$ is preserved and the proof is  complete.  
\hfill \qeds
\bigskip

\noindent\textbf{Proof of Corollary \ref{Cor:speculation}.}
The result follows readily by setting $x=y=0$ in Eq. (\ref{Eq:V_star}). 
\hfill 
\qeds


%
\noindent\textbf{Proof of Lemma \ref{Lem:increasing}.}
We prove the result by backward induction. In
particular, in each iteration, we will prove properties (1) and (2) by
recursively repeating two steps: deducing the property of $G_{n}$
from the property of $V_{n+1}$ and obtaining the property of $V_{n}$
from the property of $G_{n}$. Throughout the proof, for a matrix or a
vector $w$, we denote its transpose by $w^{T}$.
The {\it Hessian Matrix} (if exists) of a function $G=G(x,y)$ will be denoted by ${\bH}^{G}\left(x,y\right)$. For example, the Hessian Matrix of
$V_{n+1}\left(x_{n+1},y_{n+1}\right)$ is denoted by
\begin{eqnarray}\label{Eq:H-Vn+1}
{\bH}^{V_{n+1}}=\left[
\begin{array}{cc} \frac{\partial^2V_{n+1}}{\partial x_{n+1} \partial
x_{n+1} } &
\frac{\partial^2 V_{n+1}}{\partial x_{n+1} \partial y_{n+1} }\\
\frac{\partial^2 V_{n+1}}{\partial y_{n+1} \partial x_{n+1}} &
\frac{\partial^2V_{n+1}}{\partial y_{n+1} \partial y_{n+1}}\\
\end{array}\right]
\end{eqnarray}

%

\noindent 1).  For $V_N$, we have a one period problem. In this case, the
result for function $G_N(z_N,x_{N},y_{N})$ is obtained by Lemma
\ref{Lem:concave} with $z_n=x_n+q_n$ and the result for $V_N(x_{N},y_{N})$ is given by
Theorem \ref{Them:Optimal-v} of the single period problem.

\noindent 2).  For $n=1,2,\cdots,N-1$, we prove the results recursively via the following two steps:\\
{\bf Step 1}. We show that $G_{n}(z_{n}, x_{n}, y_{n})$ is
increasing in $y_{n}$ and concave in $z_n$, $x_n$ and $y_n$ if
$V_{n+1}(x_{n+1},y_{n+1})$ is increasing in $y_{n+1}$ and concave in
$x_{n+1}$ and $y_{n+1}$.

We first compute the partial
derivatives that will be used in the sequel for any given $z_n.$
By Eq. (\ref{Eq:xn}), we have:
\begin{eqnarray}
\frac{\partial x_{n+1}}{\partial z_{n}} =\frac{\partial x_{n+1}}{\partial x_{n}}&=&{\bf 1}_{\{z_n>D_n\}} \label{Eq:N_derv profit_x-qx} ,\\
 \frac{\partial x_{n+1}}{\partial y_{n}}&=&0  . \label{Eq:N_derv profit_x-y}
\end{eqnarray}
Similarly, from Eq. (\ref{Eq:yn}) we obtain:
\begin{eqnarray}
\frac{\partial y_{n+1}}{\partial z_{n}}&=&
p'_n{\bf
1}_{\{z_n<D_n\}}-h'_n{\bf 1}_{\{z_n>D_n\}}-
c'_n\left[\ (1+i_n){\bf
1}_{\{z_n<x_n+y_n\}}+
(1+\ell_n){\bf
1}_{\{z_n>x_n+y_n\}}
\right]\label{Eq:N_derv profit_y-q} ;\\
\frac{\partial y_{n+1}}{\partial x_{n}}&=&p'_n{\bf
1}_{\{z_n<D_n\}}-h'_n{\bf 1}_{\{z_n>D_n\}}\label{Eq: N_derv profit_y-x};\\
\frac{\partial y_{n+1}}{\partial y_{n}}&=&c'_n\left[(1+i_n){\bf
1}_{\{z_n<x_n+y_n\}}+(1+\ell_n){\bf 1}_{\{z_n>x_n+y_n\}}\right]\label{Eq:N_derv profit_y-y}.
\end{eqnarray}

From Eqs. (\ref{Eq:N_derv profit_x-qx})-(\ref{Eq:N_derv profit_y-y}), it readily follows that the second
order derivatives of $x_{n+1}$ and $y_{n+1}$
with respect to $z_n$, $x_n$ and $y_n$ are all zero.
In the sequel we interchange  differentiation and
integration in several places, which can be justified by the {\it Lebesgue's Dominated
Convergence} Theorem; cf. \cite{bartle1995elements}.

 The increasing property of function $G_n(z_n,x_n,y_n)$ in  $y_n$ can
be established by taking the first order derivative of  Eq.
(\ref{Eq: E profit}) with respect to $y_n$. Then, we have
\begin{eqnarray}\label{Eq:N_ derv profit_y}
\frac{\partial}{\partial y_{n}}G_{n}(z_{n},x_{n},y_{n}) &= &{\bE}\left[\frac{\partial V_{n+1}(x_{n+1},y_{n+1})}{\partial x_{n+1}}
\frac{\partial x_{n+1}}{\partial y_{n}}+ \frac{\partial
V_{n+1}(x_{n+1},y_{n+1})}{\partial y_{n+1}}
\frac{\partial y_{n+1}}{\partial y_{n}}\right]\nonumber\\
 &= &{\bE}\left[\frac{\partial V_{n+1}(x_{n+1},y_{n+1})}{\partial y_{n+1}}
\frac{\partial y_{n+1}}{\partial y_{n}}\right]\geq 0\nonumber ,
\end{eqnarray}
where the second equality holds since
$\partial x_{n+1} /\partial y_n = 0,$  by Eq. (\ref{Eq:N_derv profit_x-y}), and the inequality holds by Eq. (\ref{Eq:N_derv profit_y-y}) and
the induction hypothesis that
$V_{n+1}$ is increasing in $y_{n+1}$. 

To prove the concavity of $G_n(z_n,x_n,y_n)$ in $z_n$, $x_n$ and $y_n$, we
compute its Hessian matrix and show that it is negative semi-definite. To this end, we compute the first order
partial derivatives  of $V_{n+1}$ as,
\begin{eqnarray}
\frac{\partial G_n}{\partial z_{n}} &=& {\bE}\left[\frac{\partial V_{n+1}(x_{n+1},y_{n+1})}{\partial x_{n+1}}
\frac{\partial x_{n+1}}{\partial z_{n}}+ \frac{\partial
V_{n+1}(x_{n+1},y_{n+1})}{\partial y_{n+1}} \frac{\partial
y_{n+1}}{\partial z_{n}}\right]\label{Eq: N_derv profit_q};\\
\frac{\partial G_{n}}{\partial x_{n}} &=& {\bE}\left[\frac{\partial
V_{n+1}(x_{n+1},y_{n+1})}{\partial x_{n+1}} \frac{\partial
x_{n+1}}{\partial x_{n}}+ \frac{\partial
V_{n+1}(x_{n+1},y_{n+1})}{\partial y_{n+1}} \frac{\partial
y_{n+1}}{\partial x_{n}}\right]\label{Eq: N_V_derv1 x};\\
\frac{\partial G_{n}}{\partial y_{n}} &=& {\bE}\left[\frac{\partial
V_{n+1}(x_{n+1},y_{n+1})}{\partial x_{n+1}} \frac{\partial
x_{n+1}}{\partial y_{n}}+ \frac{\partial
V_{n+1}(x_{n+1},y_{n+1})}{\partial y_{n+1}} \frac{\partial
y_{n+1}}{\partial y_{n}}\right]\label{Eq: N_V_derv1 y},
\end{eqnarray}
and consequently the second order derivatives as
\begin{eqnarray*}
\frac{\partial ^2 G_{n}}{\partial z_{n} \partial z_{n}} &=&\bE\left[\left[ \frac{\partial x_{n+1}}{\partial
z_{n}},\frac{\partial y_{n+1}}{\partial z_{n}} \right]\cdot {\bH}^{V_{n+1}} \cdot \left[ \frac{\partial x_{n+1}}{\partial
z_{n}},\frac{\partial
y_{n+1}}{\partial z_{n}} \right]^{T}\right];\\
\frac{\partial ^2 G_{n}}{\partial z_{n} \partial x_{n}} &=&\bE\left[ \left[ \frac{\partial x_{n+1}}{\partial
z_{n}},\frac{\partial y_{n+1}}{\partial z_{n}} \right]\cdot {\bH}^{V_{n+1}} \cdot \left[ \frac{\partial x_{n+1}}{\partial
x_{n}},\frac{\partial
y_{n+1}}{\partial x_{n}} \right]^{T} \right];\\
\frac{\partial ^2 G_{n}}{\partial z_{n} \partial y_{n}} &=&\bE\left[ \left[ \frac{\partial x_{n+1}}{\partial
z_{n}},\frac{\partial y_{n+1}}{\partial z_{n}} \right]\cdot {\bH}^{V_{n+1}} \cdot \left[ \frac{\partial x_{n+1}}{\partial
y_{n}},\frac{\partial
y_{n+1}}{\partial y_{n}} \right]^{T}\right];\\
\frac{\partial ^2 G_{n}}{\partial x_{n} \partial x_{n}} &=&\bE\left[ \left[ \frac{\partial x_{n+1}}{\partial
x_{n}},\frac{\partial y_{n+1}}{\partial x_{n}} \right]\cdot {\bH}^{V_{n+1}} \cdot \left[ \frac{\partial x_{n+1}}{\partial
x_{n}},\frac{\partial
y_{n+1}}{\partial x_{n}} \right]^{T}\right];\\
\frac{\partial ^2 G_{n}}{\partial x_{n} \partial y_{n}} &=& \bE\left[ \left[ \frac{\partial x_{n+1}}{\partial
x_{n}},\frac{\partial y_{n+1}}{\partial x_{n}} \right]\cdot {\bH}^{V_{n+1}} \cdot \left[ \frac{\partial x_{n+1}}{\partial
y_{n}},\frac{\partial y_{n+1}}{\partial y_{n}} \right]^{T}\right];\\
\frac{\partial ^2 G_{n}}{\partial y_{n} \partial y_{n}} &=&\bE\left[  \left[ \frac{\partial x_{n+1}}{\partial
y_{n}},\frac{\partial y_{n+1}}{\partial y_{n}} \right]\cdot {\bH}^{V_{n+1}} \cdot \left[ \frac{\partial x_{n+1}}{\partial
y_{n}},\frac{\partial y_{n+1}}{\partial y_{n}} \right]^{T} \right],
\end{eqnarray*}
\noindent where, by Eqs. (\ref{Eq:N_derv profit_x-qx}) - (\ref{Eq:N_derv profit_y-y}), the terms involved with the second order
derivatives of $x_{n+1}$ and $y_{n+1}$ with respect to $x_{n}$ and
$y_{n}$ have vanished  and where the Hessian matrix $\bH^{V_{n+1}}$ is given by Eq. (\ref{Eq:H-Vn+1}).
Thus, the Hessian matrix of $G_n$ is:
\begin{eqnarray}\label{H-G}
{\bH}^{G_n}\left(z_n, x_n,y_n\right)=\left[
\begin{array}{ccc}
\frac{\partial^2 G_{n}}{\partial z_{n} \partial
z_{n} }& \frac{\partial^2 G_{n}}{\partial z_{n} \partial
x_{n} } & \frac{\partial^2 G_{n}}{\partial z_{n} \partial y_{n} }\\
\frac{\partial^2 G_{n}}{\partial x_{n} \partial
z_{n} }&\frac{\partial^2 G_{n}}{\partial x_{n} \partial
x_{n} } & \frac{\partial^2 G_{n}}{\partial x_{n} \partial y_{n} }\\
\frac{\partial^2 G_{n}}{\partial y_{n} \partial
z_{n} }&\frac{\partial^2 G_{n}}{\partial y_{n} \partial x_{n}} &
\frac{\partial^2 G_{n}}{\partial y_{n} \partial y_{n}}\\
\end{array}\right] ,
\end{eqnarray}
with its elements as given above. To prove it is negative
semi-definite,
we consider the quadratic
function below for any real $w_1$, $w_2$ and $w_3$,
\begin{eqnarray}\label{H-G_quadratic}
\left[w_1, w_2, w_3\right]\cdot {\bH}^{G_n}\cdot \left[w_1, w_2, w_3\right]^T&=&\bE\left[ \vec{v} \cdot \bH^{V_{n+1}} \cdot \vec{v}^T\right],
\end{eqnarray}
where $\vec{v}=w_1\left[ \frac{\partial x_{n+1}}{\partial
z_{n}},\frac{\partial y_{n+1}}{\partial z_{n}} \right]+w_2\left[ \frac{\partial x_{n+1}}{\partial
x_{n}},\frac{\partial y_{n+1}}{\partial x_{n}} \right]+w_3\left[ \frac{\partial x_{n+1}}{\partial
y_{n}},\frac{\partial y_{n+1}}{\partial y_{n}} \right]$.
Since by the induction hypothesis  ${\bH}^{V_{n+1}}$ is negative
semi-definite, we have for any $w_1$, $w_2$ and $w_3$
$$\vec{v} \cdot {\bH}^{V_{n+1}} \cdot  \vec{v}^{T }\leq 0$$
and this implies that the left side of Eq. (\ref{H-G_quadratic}) is non-positive. Thus, the proof for Step 1 is
complete.\\

{\bf Step 2}. We show that $V_{n}(x_{n}, y_{n})$ is concave in
$(x_{n}, y_{n})$ if $G_{n}(z_n,x_{n},y_{n})$ is concave.

Since $G_{n}(z_n,x_{n},y_{n})$ is concave in $z_{n}$ and $x_{n}$ and
$y_{n}$, $V_{n}(x_{n}, y_{n}) =\max_{z_n\geq x_n}
G_{n}(z_n,x_{n},y_{n})$ is concave in $x_n,$ $y_n$ by the fact that concavity is reserved under maximization; cf. Proposition A.3.10 in  \cite{zipkin2000foundations}.

Finally, the induction proof is complete with steps 1 and 2. 
\hfill \qedsymbol\\

\noindent\textbf{Proof of Theorem \ref{Them:Dynamic_order_policy}.}
Given state $(x_n,y_n)$ at the beginning of period
$n\in \bN$, we consider,
\begin{equation}\label{Eq:DerG=0}
\frac{\partial}{\partial z_{n}}G_{n}(z_{n},x_{n},y_{n}) =0,
\end{equation}
 where ${\partial G_{n}(z_{n},x_{n},y_{n})}/{\partial z_{n}} $ is given by
Eq. (\ref{Eq: N_derv profit_q}). Substituting Eqs. (\ref{Eq:N_derv
profit_x-qx}) and (\ref{Eq:N_derv profit_y-q}) into Eq. (\ref{Eq:
N_derv profit_q}) we consider the following cases:

\indent {\bf (1)} for $z_n\leq x_n+y_n$,
\begin{eqnarray}\label{Eq: N_derv profit_q2i}
\frac{\partial G_{n}}{\partial z_{n}}&=& {\bE}\left[\frac{\partial
V_{n+1}}{\partial x_{n+1}} {\bf 1}_{\{z_n>D_n\}}+
\frac{\partial V_{n+1}}{\partial y_{n+1}}\left( p'_n{\bf
1}_{\{z_n<D_n\}}-h'_n{\bf 1}_{\{z_n>D_n\}}-
c'_n(1+i_n)\right)\big| x_n,y_n\right];
\end{eqnarray}

\indent {\bf (2)} for  $z_n>x_n+y_n$,
\begin{eqnarray}\label{Eq: N_derv profit_q2l}
\frac{\partial G_{n}}{\partial z_{n}}&=& {\bE}\left[\frac{\partial
V_{n+1}}{\partial x_{n+1}} {\bf 1}_{\{z_n>D_n\}}+
\frac{\partial V_{n+1}}{\partial y_{n+1}}\left( p'_n{\bf
1}_{\{z_n<D_n\}}-h'_n{\bf 1}_{\{z_n>D_n\}}-
c'_n(1+\ell_n)\right)\big| x_n,y_n \right],
\end{eqnarray}

where for each case above, the random variables $x_{n+1}$ and $y_{n+1}$ within the expectations are given by Eqs. (\ref{Eq:xn}) and (\ref{Eq:yn}), respectively.

The results follow readily by
setting
the right sides of Eqs. (\ref{Eq: N_derv profit_q2i}) and (\ref{Eq:
N_derv profit_q2l}) equal to zero and simple manipulations.
Note that ${\partial G_{n}(z_{n},x_{n},y_{n})}/{\partial z_{n}} $ is monotonically decreasing in $z_n$ due to its concavity, shown in part (1) of Lemma \ref{Lem:increasing}. Consequently, there are unique solutions to each of these equations. 
\hfill \qedsymbol\\
\bigskip

\noindent\textbf{Proof Corollary \ref{Cor:alpha-beta-x+y}.}
For period $N$, the independence of $x_N$ or $y_N$ is obvious since this is a single period. For period $n<N$, revisit
Eqs. (\ref{Eq: N_derv profit_q2i}) and (\ref{Eq: N_derv profit_q2l}).   Note that $x_{n+1}$ is independent of $(x_n,y_n)$ by Eq. (\ref{Eq:xn}) while $y_{n+1}$ is dependent of $x_n+y_n$ by Eq. (\ref{Eq:yn}). Therefore, $\alpha_n$ and $\beta_n$ implicitly given by Eqs. (\ref{Eq: N_derv profit_q2i}) and (\ref{Eq: N_derv profit_q2l}) are dependent on $\xi_n=x_n+y_n$ only, and thus completes the proof.
\hfill \qedsymbol\\
\bigskip

\noindent\noindent \textbf{Proof of Lemma \ref{Lem:E-inequality}.}
To prove (a), we only give the proof for the case of increasing $f(x)$ and $g(x)$. The same argument can be applied for the case of decreasing $f(x)$ and $g(x)$.\\
Let $X'$ be another random variable which is i.i.d. of $X$. Since $f(x)$ and $g(x)$ are increasing, we always have
$[f(X)-f(X')][g(X)-g(X')]\geq 0.$
Taking expectations with respect to $X$ and $X'$ yields
\begin{eqnarray}
&&\bE\big[[f(X)-f(X')][g(X)-g(X')]\big] \nonumber\\
&=&\bE[f(X)g(X)+f(X')g(X')-f(X')g(X)-f(X)g(X')]\nonumber\\
                             &=&\bE[f(X)g(X)]+\bE[f(X')g(X')]-\bE[f(X')]\bE[g(X)]-\bE[f(X)]\bE[g(X')]\nonumber\\
                             &=&2\bE[f(X)g(X)]-2\bE[f(X)]\bE[g(X)]\geq 0.\nonumber
\end{eqnarray}
The result of part (a) readily follows from the above. \\
In a similar vein, we can prove part (b) via changing the direction of the inequality above. 
\hfill \qedsymbol\\
\bigskip

\noindent\textbf{Proof of Theorem \ref{Thm:Optimal-myopical-hat}.}
We only prove the result for $\alpha_n$. The same argument can be applied to prove the result for $\beta_n$.\\
\noindent In view of Eq. (\ref{Eq: N_derv profit_q2i}), $\alpha_n$ is the unique solution of the equation below,
\begin{eqnarray}\label{Eq:alpha_n_a}
{\bE}\left[\frac{\partial V_{n+1}}{\partial y_{n+1}}\left( p'_n{\bf
1}_{\{\alpha_n<D_n\}}-h'_n{\bf 1}_{\{\alpha_n>D_n\}}-
c'_n(1+\ell_n)\right)\right]=-{\bE}\left[\frac{\partial
V_{n+1}}{\partial x_{n+1}} {\bf 1}_{\{\alpha_n>D_n\}}\right].
\end{eqnarray}

Since $\frac{\partial V_{n+1}}{\partial x_{n+1}}\geq 0$ by Lemma \ref{Lem:increasing} part (2), the equation above is negative, which implies
\begin{eqnarray}\label{Eq:alpha_n_b}
{\bE}\left[\frac{\partial V_{n+1}}{\partial y_{n+1}}\left( p'_n{\bf
1}_{\{\alpha_n<D_n\}}-h'_n{\bf 1}_{\{\alpha_n>D_n\}}-
c'_n(1+\ell_n)\right)\right] \leq 0.
\end{eqnarray}

\noindent Further note that for any realization of demand $D_{n}=d>0$, the two terms  of the  left hand side of Eq. (\ref{Eq:alpha_n_b}):
$$\frac{\partial V_{n+1}(x_{n+1}(d),y_{n+1}(d))}{\partial y_{n+1}(d)}$$
and
$  p'_n{\bf
1}_{\{\alpha_n<d\}}-h'_n{\bf 1}_{\{\alpha_n>d\}}-
c'_n(1+\ell_n) $
are both increasing in $d$. Specifically, the first term is increasing by the concavity of $V_{n+1}$ [cf. Lemma \ref{Lem:increasing}, part (2)] and Eq. (\ref{Eq:yn}).  Then, by Lemma \ref{Lem:E-inequality} and Eq. (\ref{Eq:alpha_n_b}), one has,
\begin{eqnarray}\label{Eq:alpha_n_c}
{\bE}\left[\frac{\partial V_{n+1}}{\partial y_{n+1}}\right]{\bE}\left[ p'_n{\bf
1}_{\{\alpha_n<D_n\}}-h'_n{\bf 1}_{\{\alpha_n>D_n\}}-
c'_n(1+\ell_n)\right]\leq 0
\end{eqnarray}

\noindent Since $\frac{\partial V_{n+1}}{\partial y_{n+1}}\geq 0$    by Lemma \ref{Lem:increasing} part (2), 
the above inequality implies
\begin{eqnarray}\nonumber
{\bE}\left[ p_n{\bf
1}_{\{\alpha_n<D_n\}}-h_n{\bf 1}_{\{\alpha_n>D_n\}}-
c_n(1+\ell_n)\right]\leq 0,
\end{eqnarray}
which, after simple algebra,  is equivalent to
$p_n-c_n\cdot (1+\ell_n)-(p_n-s_n)F_n(\alpha_n)\leq 0.$
The above further simplifies to
$$F(\alpha_n)\geq \frac{p_n-c_n\cdot (1+\ell_n)}{p_n-s_n}.$$
By Eqs. (\ref{Eq:alpha_hat}) and (\ref{Eq:myoptical-alpha-hat}), the right hand side in the above inequality is $F_n({\underline \alpha}_n)$. Thus, we have $F_n(\alpha_n)\geq F_n({\underline \alpha}_n)$, which completes the proof for $\alpha_n\geq {\underline \alpha}_n$ by the increasing property of $F_n(\cdot)$.  
\hfill \qedsymbol\\
%
\bigskip

\noindent\textbf{Proof of Proposition \ref{Pro:x-y}.}
To prove part (i), it suffices to prove that $V_n(x_n,y_n)\leq V_n(x_n-d,y_n+d)$ for arbitrarily small $d>0$. To this end, consider the initial state $(x_n-d,y_n+d)$. Then, the \NV\ can always purchase $d$ units without any additional cost to reset the initial state to be $(x_n,y_n)$. This means $V_n(x_n,y_n)\leq V_n(x_n-d,y_n+d)$.\\
The proof for part (ii) follows from part (i). Considering the derivative of $V_n(x_n-d,y_n+d)$ with respect to $d$, we have
\begin{eqnarray}\nonumber
\lim_{d\to 0}\frac{V_n(x-d,y+d)-V_n(x,y)}{d}=-\partial V_n(x,y)/\partial x+\partial V_n(x,y)/\partial y\geq 0
\end{eqnarray}
where the inequality follows from part (i). Finally, part (ii) readily follows by rearranging the above equation.
\hfill \qedsymbol\\
\bigskip

\noindent\textbf{Proof of Proposition \ref{Pro:aL}}.
To prove Part (i), we note that:
 $$
V^{L}_n(x_{n},y_{n})
=
 \max_{z_n\geq x_n}{\bE}\left[ V_{n+1}(x_{n+1}, y_{n+1} ) \right]
\leq {\bE}\left[ V_{n+1}(0,x_{n+1}+y_{n+1})\right]
=
V^{L}_n(x_{n},y_{n})
$$
where the inequality above follows from Proposition \ref{Pro:x-y}(i) for period  $n+1 $ with $d= x_{n+1} $.
 
For part (ii), we only show that  $\alpha^L_n\geq \alpha_n$ since the same argument can be applied to prove $\beta^L_n\geq \beta_n$.

To begin with, note that by Lemma \ref{Lem:increasing} one can show that 
    $\bE[V_{n+1}(0,\xi_{n+1}(z_n))]$ is concave in $z_n$ and its derivative is decreasing in $z_n$. Therefore, 
     $z_n=\alpha^L_n$ is the unique solution of
    $\bE[V_{n+1}(0,\xi_{n+1}(z_n))]=0$. In addition,  the derivative of    $\bE[V_{n+1}(0,\xi_{n+1}(z_n))]$ is  positive for $z_n< \alpha^L_n$ and negative for $z_n> \alpha^L_n$. Therefore, it suffices to show that at $z_n=\alpha_n$ the following is true
\begin{eqnarray}\label{Eq:dVz_alpha}
\frac{d}{dz_n} \bE\bigg[V_{n+1}\big(0,\, x_{n+1}(z_n)+y_{n+1}(z_n)\big)\bigg]\bigg |_{z_n=\alpha_n}\geq 0.
\end{eqnarray}
 It is easy to show by backward induction that,  
 $d V_{n+1}\big(0,\, x_{n+1}(z_n)+y_{n+1}(z_n)\big)/dz_{n}$ is bounded for   every $z_{n}$ in terms of a function of $D_{n}$, (the last period $R_{N}(D_N,q_{N},x_{N})+K_{N}(q_{N},y_{N})$ is bounded for any $D_N$).
The previous condition allows the interchange between differentiation and integration; cf. \cite{Widder1989}. Accordingly, we have the following,
\begin{eqnarray}\label{Eq: N_derv profit_bound}
\frac{d \bE [V_{n+1}(0,\xi_{n+1})]}{d z_{n}}&=&\bE [\frac{d V_{n+1}(0,\xi_{n+1})}{d z_{n}}]\nonumber\\
&=& {\bE}\left[\frac{\partial
V_{n+1}}{\partial \xi_{n+1}} {\bf 1}_{\{z_n>D_n\}}+
\frac{\partial V_{n+1}}{\partial \xi_{n+1}}\left( p'_n{\bf
1}_{\{z_n<D_n\}}-h'_n{\bf 1}_{\{z_n>D_n\}}-
c'_n(1+\ell_n)\right)\right],\nonumber\\
&\geq& {\bE}\left[\frac{\partial
V_{n+1}}{\partial x_{n+1}} {\bf 1}_{\{z_n>D_n\}}+
\frac{\partial V_{n+1}}{\partial y_{n+1}}\left( p'_n{\bf
1}_{\{z_n<D_n\}}-h'_n{\bf 1}_{\{z_n>D_n\}}-
c'_n(1+\ell_n)\right) \right]=0,\nonumber
\end{eqnarray}
where for simplicity we write $\frac{\partial
V_{n+1}}{\partial \xi_{n+1}}:=\frac{\partial
V_{n+1}(0,\,\xi_{n+1})}{\partial \xi_{n+1}}$, 
$\frac{\partial
V_{n+1}}{\partial x_{n+1}}:=\frac{\partial
V_{n+1}(x_{n+1},\,y_{n+1})}{\partial x_{n+1}}$ and 
$\frac{\partial
V_{n+1}}{\partial y_{n+1}}:=\frac{\partial
V_{n+1}(x_{n+1},\,y_{n+1})}{\partial y_{n+1}}$.

The inequality above is justified using backward induction as follows. 

\noindent First, since  the marginal contribution of profit of an additional unit of capital for $V_{n+1}(0,\xi_{n+1})$ is larger than that of $V_{n+1}(x_{n+1},y_{n+1})$ where $\xi_{n+1}=x_{n+1}+y_{n+1}$, we have:
\begin{eqnarray}
\partial V_{n+1}(x_{n+1},\,y_{n+1})/\partial y_{n+1}&\leq &\partial V_{n+1}(0,\,\xi_{n+1})/\partial \xi_{n+1}; \nonumber\\
\partial V_{n+1}(x_{n+1},\,y_{n+1})/\partial x_{n+1}&\leq &\partial V_{n+1}(0,\,\xi_{n+1})/\partial \xi_{n+1}.\nonumber
\end{eqnarray}
Second, since  $z_n=\alpha_n$ is the unique optimum of $\bE[V_{n+1}(x_{n+1},y_{n+1})]$, its derivative at $z_n=\alpha_n$ is zero. This proves Eq. (\ref{Eq:dVz_alpha}) and thus completes the proof.
\hfill \qedsymbol\\

\bigskip

\noindent\textbf{Proof of Proposition \ref{Pro:Optimal-myopical-tilde}.}
For period $N$, the result readily follows from the optimal solution of the single-period model.  For all other periods, we only prove the result for ${\bar \alpha}_n\geq \alpha_n$ since a similar argument can be used to prove ${\bar \beta}_n\geq \beta_n$ by replacing $\ell_n$ with $i_n$.


\noindent By Proposition \ref{Pro:aL}, we have $\alpha_n\leq \alpha^{L}_{n}$ and $\alpha^{L}_{n}$   is by differentiating Eq. (\ref{Eq:tildeV-e}) and setting the derivative to zero, that is,
\begin{eqnarray}\label{Eq:V(x+y)}
{\bE}\left[
\frac{\partial
V_{n+1}(0,\xi_{n+1})}{\partial \xi_{n+1}}\left( {\bf 1}_{\{{\alpha^{L}_{n}}>D_n\}}+ p'_n{\bf
1}_{\{{ \alpha^{L}_n}<D_n\}}-h'_n{\bf 1}_{\{{ \alpha^{L}_n}>D_n\}}-
c'_n(1+\ell_n)\right)\right]=0.
\end{eqnarray}
For any realization of demand $D_{n}=d>0$ the term
${\partial V_{n+1}(0,\xi_{n+1}(d)\,)}/{\partial \xi_{n+1}(d)}$ is decreasing in $d$     by the concavity of $V_{n+1}$ [cf. Lemma \ref{Lem:increasing} part (2)] and the fact that $\xi_{n+1}$ is increasing in $d$ by Eqs. (\ref{Eq:xn})-(\ref{Eq:yn}).

Furthermore, the term
$ {\bf 1}_{\{{ \alpha}_n>d\}}+ p'_n{\bf
1}_{\{{ \alpha}_n<d\}}-h'_n{\bf 1}_{\{{ \alpha}_n>d\}}-
c'_n(1+\ell_n)=  p'_n-(p'_n+h'_n-1){\bf 1}_{\{{ \alpha}_n>d\}}-
c'_n(1+\ell_n),$
is increasing in $d$.
By Eq. (\ref{Eq:V(x+y)}) and Lemma \ref{Lem:E-inequality}, one has
\begin{eqnarray}\nonumber
{\bE}\left[
\frac{\partial
V_{n+1}(0,\xi_{n+1})}{\partial \xi_{n+1}}\right]\cdot {\bE}\left[ {\bf 1}_{\{{\alpha^{L}_{n}}>D_n\}}+ p'_n{\bf
1}_{\{{ \alpha^{L}_n}<D_n\}}-h'_n{\bf 1}_{\{{ \alpha^{L}_n}>D_n\}}-
c'_n(1+\ell_n)\right]\geq 0.
\end{eqnarray}
 Since $
{\partial
V_{n+1}(0,\xi_{n+1})}/{\partial \xi_{n+1}}\geq 0$  by Lemma \ref{Lem:increasing} part (2),
the above inequality implies
\begin{eqnarray}\label{Eq:alpha_x+y_b}
 {\bE}\left[ {\bf 1}_{\{{\alpha^{L}_{n}}>D_n\}}+ p'_n{\bf
1}_{\{{ \alpha^{L}_n}<D_n\}}-h'_n{\bf 1}_{\{{ \alpha^{L}_n}>D_n\}}-
c'_n(1+\ell_n)\right]\geq 0,
\end{eqnarray}
which, after simple algebra,  is equivalent to
$p_n-c_n\cdot (1+\ell_n)-(p_n+h_n-c_{n+1})F_n(\alpha^L_n)\geq 0.$
The above further simplifies to
$$F(\alpha^L_n)\leq \frac{p_n-c_n\cdot (1+\ell_n)}{p_n+h_n-c_{n+1}}.$$
Note that the right hand side of the above is less than 1 since $c_n(1+\ell_n)+h_n\geq c_{n+1}$ by assumption. Next, by Eqs. (\ref{Eq:tilde-a}) and (\ref{Eq:myoptical-alpha}), the right hand side in the above inequality is $F_n({\bar \alpha}_n)$. Thus, we have $ F_n({\bar \alpha}_n) \geq F_n(\alpha^L_n)$, which means ${ \bar \alpha}_n \geq \alpha^L_n $. Thus, the proof for $\bar \alpha_n\geq { \alpha}_n$ is complete, since  $\alpha^L_n\geq {  \alpha}_n$ by Proposition \ref{Pro:aL}.  
\hfill \qedsymbol\\
\bigskip

\textbf{Proof of Theorem \ref{Them:SellingBack}.}
Since $V^S_n(\xi)$ is increasing and concave in $\xi$,
to prove part (i), we can follow the proof for Theorem \ref{Them:Dynamic_order_policy} by simplifying the two-variable state $(x_n,y_n)$ as $\xi_n$.  The two threshold values $\alpha^S_n$ and $\beta^S_n$ can be obtained via solving the first order conditions under loan with rate $\ell_n$ and deposit with rate $r_n$, respectively.

To prove part (ii), we only need to prove the first inequality since the second one has been proved by Proposition \ref{Pro:aL}. This can be done via backward induction. First, for the induction basis, note that for the last period $N$,  one has $V^S_n(x+y)=V^L_n(x,y)$. 
Next, for the induction step, assume that the inequality holds for periods $n+1,n+2, \ldots, N$. Then, for any  $\xi_{n+1}=x_{n+1}+y_{n+1}$, one has
$V^S_{n+1}(\xi_{n+1})\geq V^L_{n+1}(x_{n+1},y_{n+1})$. 
Noting that $\xi_{n+1}$ is a function of $\xi_n$, $z_n$ and $D_n$. Taking expectations w.r.t. $D_n$ yields   
$\bE[V^S_{n+1}(\xi_{n+1})]\geq \bE[V^L_{n+1}(x_{n+1},y_{n+1})]$. Finally, taking the maximum over $z_n$, we obtain that 
the inequality  
$V^S_{n}(\xi_{n}) \geq V^L_{n}(x_{n},y_{n})$
 holds, and  the proof is complete.
\hfill \qedsymbol\\

\bigskip

\textbf{Proof of Theorem \ref{Corr:Optimal-q-piecewise}.}
The proof readily follows via a straightforward modification of the proof of Theorem \ref{Them:Optimal-q}.
\hfill \qedsymbol\\

\textbf{Proof of Theorem \ref{Them:Dynamic_order_policy_Loan_Capacity}.}
The proof follows the same steps as in Theorem \ref{Them:Optimal-q} and Theorem \ref{Them:Dynamic_order_policy} counterparts, with some minor modifications.
\hfill \qedsymbol\\
\bigskip

\textbf{Proof of Theorem \ref{Them:Dynamic_order_policy-backorder}.}
The proof follows readily from the proof for Theorem \ref{Them:Dynamic_order_policy} by replacing the inventory-cash flow equations to be those of Eqs. (\ref{Eq:xn0-backorder})-(\ref{Eq:yn0-backorder}).
\hfill \qedsymbol\\

\end{document}